\documentclass[12pt]{article}
\usepackage{amsmath,amsfonts,amssymb,amsthm}
\input amssym.def
\begin{document}
\def \Z{\Bbb Z}
\def \C{\Bbb C}
\def \R{\Bbb R}
\def \Q{\Bbb Q}
\def \N{\Bbb N}

\def \A{{\mathcal{A}}}
\def \D{{\mathcal{D}}}
\def \E{{\mathcal{E}}}
\def \E{{\mathcal{E}}}
\def \H{\mathcal{H}}
\def \S{{\mathcal{S}}}
\def \wt{{\rm wt}}
\def \tr{{\rm tr}}
\def \span{{\rm span}}
\def \Res{{\rm Res}}
\def \Der{{\rm Der}}
\def \End{{\rm End}}
\def \Ind {{\rm Ind}}
\def \Irr {{\rm Irr}}
\def \Aut{{\rm Aut}}
\def \GL{{\rm GL}}
\def \Hom{{\rm Hom}}
\def \mod{{\rm mod}}
\def \ann{{\rm Ann}}
\def \ad{{\rm ad}}
\def \rank{{\rm rank}\;}
\def \<{\langle}
\def \>{\rangle}

\def \g{{\frak{g}}}
\def \h{{\hbar}}
\def \k{{\frak{k}}}
\def \sl{{\frak{sl}}}
\def \gl{{\frak{gl}}}

\def \be{\begin{equation}\label}
\def \ee{\end{equation}}
\def \bex{\begin{example}\label}
\def \eex{\end{example}}
\def \bl{\begin{lem}\label}
\def \el{\end{lem}}
\def \bt{\begin{thm}\label}
\def \et{\end{thm}}
\def \bp{\begin{prop}\label}
\def \ep{\end{prop}}
\def \br{\begin{rem}\label}
\def \er{\end{rem}}
\def \bc{\begin{coro}\label}
\def \ec{\end{coro}}
\def \bd{\begin{de}\label}
\def \ed{\end{de}}

\newcommand{\m}{\bf m}
\newcommand{\n}{\bf n}
\newcommand{\nno}{\nonumber}
\newcommand{\nord}{\mbox{\scriptsize ${\circ\atop\circ}$}}
\newtheorem{thm}{Theorem}[section]
\newtheorem{prop}[thm]{Proposition}
\newtheorem{coro}[thm]{Corollary}
\newtheorem{conj}[thm]{Conjecture}
\newtheorem{example}[thm]{Example}
\newtheorem{lem}[thm]{Lemma}
\newtheorem{rem}[thm]{Remark}
\newtheorem{de}[thm]{Definition}
\newtheorem{hy}[thm]{Hypothesis}
\makeatletter
\@addtoreset{equation}{section}
\def\theequation{\thesection.\arabic{equation}}
\makeatother
\makeatletter

\begin{center}
{\Large \bf On vertex algebras and their modules associated with
even lattices}
\end{center}

\begin{center}
{Haisheng Li\footnote{Partially supported
by NSF grant DMS-0600189} and Qing Wang\footnote{Permanent address:
Department of Mathematics, Xiamen University, Xiamen, China}\\
Department of Mathematical Sciences, Rutgers University, Camden, NJ
08102}
\end{center}

\begin{abstract}
We study vertex algebras and their modules associated with possibly
degenerate even lattices, using an approach somewhat different from
others. Several known results are recovered and a number of new
results are obtained. We also study modules for Heisenberg algebras
and we classify irreducible modules satisfying certain conditions
and obtain a complete reducibility theorem.
\end{abstract}

\section{Introduction}
Let $L$ be any nondegenerate even lattice in the sense that $L$ is a
finite-rank free abelian group equipped with a symmetric
$\Z$-bilinear form $\<\cdot,\cdot\>$ such that $\<\alpha,\alpha\>\in
2\Z$ for $\alpha\in L$. Associated to $L$ there exists a canonical
vertex algebra $V_{L}$ (\cite{b86}, \cite{flm}). These vertex
algebras were originated from the explicit vertex-operator
realizations of the basic modules for affine Kac-Moody Lie algebras
and they form an important class in the theory of vertex algebras.
For vertex algebras $V_{L}$, irreducible modules were classified in
\cite{dong1} (cf. \cite{dlm-reg}, \cite{ll}) where Dong proved that
every irreducible $V_{L}$-module is isomorphic to one of those
constructed in \cite{flm}.

The vertex algebras $V_{L}$ were built up by using certain
infinite-dimensional Heisenberg (Lie) algebras and their modules.
For a nondegenerate even lattice $L$, set $\frak{h}=\C\otimes_{\Z}L$
and extend the $\Z$-bilinear form on $L$ to a nondegenerate
symmetric $\C$-bilinear form $\<\cdot,\cdot\>$ on $\frak{h}$.
Associated to the pair $(\frak{h},\<\cdot,\cdot\>)$, there is an
affine Lie algebra $\hat{\frak{h}}=\frak{h}\otimes
\C[t,t^{-1}]\oplus \C {\bf k}$, whose subalgebra
$\hat{\frak{h}}_{*}=\sum_{n\ne 0}\frak{h}(n)+\C {\bf k}$ is a
Heisenberg algebra, where $\frak{h}(m)=\frak{h}\otimes t^{m}$ for
$m\in \Z$. For each $\alpha\in \frak{h}=\frak{h}^{*}$, there is a
canonical irreducible $\hat{\frak{h}}$-module $M(1,\alpha)$ on which
${\bf k}$ acts as identity and $h(0)$ acts as scalar $\<\alpha,h\>$
for $h\in \frak{h}$. The associated vertex algebra $V_{L}$ is built
on the direct sum of the non-isomorphic irreducible
$\hat{\frak{h}}$-modules $M(1,\alpha)$ for $\alpha\in L$, where
$M(1,0)$ is a vertex subalgebra with $M(1,\alpha)$ as irreducible
modules. Note that the vertex algebra $M(1,0)$ can also be
constructed independently by starting with the finite-dimensional
vector space $\frak{h}$ equipped with a nondegenerate symmetric
bilinear form $\<\cdot,\cdot\>$, without referring to the lattice
$L$. All of these motivated the work \cite{lx}.

The main goal of \cite{lx} was to characterize the family of vertex
algebras $V_{L}$. With vertex algebras $V_{L}$ as models, a class
$\mathcal{A}$ of simple vertex algebras was formulated, where each
$V$ in $\mathcal{A}$ contains $M(1,0)$ associated to some $\frak{h}$
as a vertex subalgebra and is a direct sum of some non-isomorphic
irreducible $M(1,0)$-modules $M(1,\alpha)$ with $\alpha\in
\frak{h}$. It was proved therein that each $V$ in $\mathcal{A}$ is a
direct sum of $M(1,\alpha)$ with $\alpha\in L$ for some additive
subgroup $L$ of $\frak{h}$ such that $\<\alpha,\alpha\>\in 2\Z$ for
$\alpha\in L$ and such that there exists a function (normalized
$2$-cocycle) $\varepsilon: L\times L\rightarrow \C^{\times}$
satisfying the condition that
\begin{eqnarray*}
& &\varepsilon(\alpha,0)=\varepsilon(0,\alpha)=1,\\
& &\varepsilon(\alpha,\beta+\gamma)\varepsilon(\beta,\gamma)
=\varepsilon(\alpha+\beta,\gamma)\varepsilon(\alpha,\beta),\\
& &\varepsilon(\alpha,\beta)\varepsilon(\beta,\alpha)^{-1}
=(-1)^{\<\alpha,\beta\>} \ \ \ \ \mbox{ for }\alpha,\beta\in L.
\end{eqnarray*}
Such a vertex algebra $V$ was denoted therein by $V_{(\frak{h},L)}$.
(However, the existence and uniqueness of the desired vertex algebra
$V_{(\frak{h},L)}$ was neglected.) If the vertex algebra $V$ is
finitely generated, one can readily see that $L$ is finitely
generated, so that $L$ is free of finite rank (as $L$ is
torsion-free). In this case, $L$ is a (finite rank) possibly
degenerate even lattice.

In another work \cite{bdt}, motivated by a certain connection
between vertex algebras and toroidal Lie algebras, Berman, Dong and
Tan studied modules for a class of vertex algebras which may be
identified as $V_{(\frak{h},L)}$ with $\<\cdot,\cdot\>|_{L}=0$ and
$\dim \frak{h}=2{\rm rank} (L)$ (finite). Among other results, they
constructed and classified $\N$-graded irreducible modules.

In this current paper, we study vertex algebra $V_{(\frak{h},L)}$
and its modules associated with a general pair $(\frak{h},L)$ {}from
the point of view of vertex algebra extensions. Specifically, we
make use of the vertex algebra $M(1,0)$ and its modules to construct
vertex algebra $V_{(\frak{h},L)}$ and its modules. This idea goes
back to the so-called simple current extension of vertex operator
algebras (\cite{li-phys}, \cite{li-ext}, \cite{dlm-simple},
\cite{li-extaffine}). For the vertex algebra $V_{L}$ with $L$ a
nondegenerate even lattice, irreducible modules have been
constructed before in \cite{flm} and \cite{ll} with quite different
methods. Our treatment here is further different from those of
\cite{flm} and \cite{ll}.

It is known (cf. \cite{ll}) that on every $\hat{\frak{h}}$-module
$W$ of level $1$, which is {\em restricted} in the sense that for
every $w\in W$, $\frak{h}(n)w=0$ for $n$ sufficiently large, there
exists a unique module structure for vertex algebra $M(1,0)$,
extending the action of $\hat{\frak{h}}$ in a certain canonical way
and that every module for vertex algebra $M(1,0)$ is naturally a
restricted $\hat{\frak{h}}$-module of level $1$. This leads us to
study irreducible restricted modules for Heisenberg algebras, e.g.,
$\hat{\frak{h}}_{*}$. It is well known that Heisenberg algebras have
canonical realizations on polynomial algebras by differential
(annihilation) operators and left multiplication (creation)
operators. For a Heisenberg algebra with a fixed nonzero level,
there is a very nice module category in which there is only one
irreducible module up to equivalence and each module is completely
reducible; see \cite{lw2} (\cite{flm}, Theorem 1.7.3), \cite{kac1}.
In the case with the Heisenberg algebra $\hat{\frak{h}}_{*}$, the
modules in this category are restricted and $\frak{h}(n)$ are
locally nilpotent for $n\ge 1$. In this paper, we study a bigger
category ${\mathcal{F}}$ of restricted modules on which
$\frak{h}(n)$ are locally finite for $n\ge 1$. It turns out that
${\mathcal{F}}$ is also very nice in the sense that each module is
completely reducible and each irreducible representation can be
explicitly realized by using differential operators and left
multiplication operators on a space of exponential functions and
polynomials.  In contrast to the old case, in the category
${\mathcal{F}}$ there are infinitely many unequivalent irreducible
objects. We also study and construct certain irreducible restricted
modules outside the category ${\mathcal{F}}$.

For a general pair $(\frak{h},L)$ with $\frak{h}\ne {\rm span}(L)$,
the vertex algebra $V_{(\frak{h},L)}$ like $M(1,0)$ has irreducible
modules other than those constructed {}from the canonical
realization of Heisenberg algebras on polynomial algebras. We
construct and classify irreducible $V_{(\frak{h},L)}$-modules in
terms of irreducible restricted $\hat{\frak{h}}$-modules of a
certain type. For the construction of irreducible modules, we use
the idea from \cite{li-ext} and \cite{dlm-simple} and for the
classification of irreducible modules, we use certain important
ideas (and results) from \cite{dong1} and \cite{dlm-reg}. In the
case that $\<\cdot,\cdot\>|_{L}=0$ we obtain certain irreducible
$V_{(\frak{h},L)}$-modules other than those constructed by Berman,
Dong and Tan in \cite{bdt}. Furthermore, we give another
construction of the vertex algebra $V_{(\frak{h},L)}$ by using a
certain affine Lie algebra.

This paper is organized in the following manner. In Section 2, we
study certain categories of modules for Heisenberg algebras and we
classify all the irreducible objects. In Section 3, we study the
modules for vertex algebras associated with even lattices. In
Section 4, we give a characterization of vertex algebras in terms of
an affine Lie algebra.

\section{Modules for Heisenberg Lie algebras}
In this section we first associate a Heisenberg Lie algebra $\H_{I}$
to each nonempty set $I$ and we then construct and classify
irreducible $\H_{I}$-modules of certain types. Especially, we define
a category ${\mathcal{F}}$ and we establish a complete reducibility
theorem, which generalizes a theorem of \cite{flm} (Section 1.7).

First, we start with introducing the Heisenberg Lie algebra
$\H_{I}$. Let $I$ be any nonempty set, which is fixed throughout
this section. Let $\H_{I}$ be the Heisenberg Lie algebra with a
basis $\{{\bf p}_{i},{\bf q}_{i},{\bf k}\;|\; i\in I\}$ (over $\C$)
and with Lie bracket relations
$$[{\bf k},\H_{I}]=0,\ \
[{\bf p}_{i},{\bf p}_{j}]=0=[{\bf q}_{i},{\bf q}_{j}]\ \ \mbox{ and
}\ \ [{\bf p}_{i},{\bf q}_{j}]=\delta_{ij}{\bf k}\ \ \ \mbox{ for }
i,j\in I.$$ An $\H_{I}$-module $W$ is said to be of {\em level}
$\ell\in \C$ if ${\bf k}$ acts as scalar $\ell$.

\bl{lnilp-finite} Let $W$ be any irreducible $\H_{I}$-module of
level $1$. For $i\in I$, ${\bf p}_{i}$ has an eigenvector in $W$ if
and only if ${\bf p}_{i}$ acts locally finitely. Furthermore, ${\bf
p}_{i}w=0$ for some nonzero $w\in W$ if and only if ${\bf p}_{i}$
acts locally nilpotently. The same assertions hold with ${\bf
p}_{i}$ replaced by ${\bf q}_{i}$.
 \el

\begin{proof} It is clear that if ${\bf p}_{i}$ is locally finite,
${\bf p}_{i}$ has an eigenvector, as the scalar field is $\C$. Now,
assume that ${\bf p}_{i}$ has an eigenvector $w_{0}$ of eigenvalue
$\lambda$. Set
$$A=\sum_{j\in I,\; j\ne i}(\C {\bf p}_{j}+\C {\bf q}_{j})+\C {\bf
p}_{i}+\C {\bf k},$$ a Lie subalgebra of $\H_{I}$. We have
$\H_{I}=\C {\bf q}_{i}\oplus A$. Since $W$ is an irreducible
$\H_{I}$-module, $W=U(\H_{I})w_{0}=U(A)\C[{\bf q}_{i}]w_{0}$. It
follows from induction that for any nonnegative integer $k$,
$\sum_{n=0}^{k}\C{\bf q}_{i}^{n}w_{0}$ is a finite-dimensional
subspace which is closed under the action of ${\bf p}_{i}$.
Furthermore, for any $a\in U(A)$, as $[{\bf p}_{i},A]=0$,
$\sum_{n=0}^{k}a\C {\bf q}_{i}^{n}w_{0}$ is a finite-dimensional
subspace which is closed under the action of ${\bf p}_{i}$. Now, it
follows that ${\bf p}_{i}$ is locally finite. For the nilpotent
case, it is also clear.
\end{proof}

\bl{lpiqi} Let $W$ be an irreducible $\H_{I}$-module of level $1$
and let $i\in I$. If ${\bf p}_{i}$ is locally finite on $W$, then
any nonzero ${\bf q}_{i}$-stable subspace is infinite-dimensional.
In particular, ${\bf q}_{i}$ is not locally finite. The same
assertion holds when ${\bf p}_{i}$ and ${\bf q}_{i}$ are
exchanged.
\el

\begin{proof} Suppose that there exists a nonzero finite-dimensional subspace
$U$ which is ${\bf q}_{i}$-stable. Then there exists $0\ne w\in U$
such that ${\bf q}_{i}w=\alpha w$ for some $\alpha\in \C$. Using
induction we get
$${\bf q}_{i}{\bf p}_{i}^{n}w=(\alpha {\bf p}_{i}^{n}-n{\bf p}_{i}^{n-1})w
\ \ \ \mbox{ for }n\ge 0.$$ It follows that ${\bf p}_{i}^{n}w\ne 0$
for $n\ge 0$.  As ${\bf p}_{i}$ is locally finite, $\C[{\bf
p}_{i}]w$ is finite-dimensional. Let $n=\dim \C[{\bf p}_{i}]w\ge 1$.
Then $w, {\bf p}_{i}w,\dots, {\bf p}_{i}^{n-1}w$ are linearly
independent and
$${\bf p}_{i}^{n}w=c_{n-1}{\bf
p}_{i}^{n-1}w+c_{n-2}{\bf p}_{i}^{n-2}w+\cdots + c_{0}w$$ for some
$c_{0},\dots,c_{n-1}\in \C$. Applying ${\bf q}_{i}$ we get
$$\alpha {\bf p}_{i}^{n}w-n{\bf p}_{i}^{n-1}w
=\alpha c_{n-1}{\bf p}_{i}^{n-1}w+u$$ for some $u\in
\sum_{r=0}^{n-2}\C{\bf p}_{i}^{r}w$. Then we obtain
$$n{\bf p}_{i}^{n-1}w=\alpha (c_{n-2}{\bf p}_{i}^{n-2}w+\cdots + c_{0}w) -u.$$
This is a contradiction as $w, {\bf p}_{i}w,\dots, {\bf
p}_{i}^{n-1}w$ are linearly independent.
\end{proof}

\bd{dIWpq} {\em Let $W$ be an $\H_{I}$-module of level $1$. Set
\begin{eqnarray*}
I_{W}^{\bf p}&=&\{ i\in I\;|\; {\bf p}_{i}\ \mbox{ is locally
finite}\},\\
I_{W}^{\bf q}&=&\{ i\in I\;|\; {\bf q}_{i}\ \mbox{ is locally
finite}\}.
\end{eqnarray*}} \ed

For our study in the next section on modules for vertex algebras
associated to Heisenberg Lie algebras we are interested in
$\H_{I}$-modules $W$ such that for every $w\in W$, ${\bf p}_{i}w=0$
for all but finitely many $i\in I$. For convenience, we call such an
$\H_{I}$-module a {\em restricted} module. In view of Lemmas
\ref{lnilp-finite} and \ref{lpiqi}, if $W$ is an irreducible
restricted $\H_{I}$-module, then $I_{W}^{\bf p}$ is a cofinite
subset of $I$ and $I_{W}^{\bf p}\cap I_{W}^{\bf q}=\emptyset$.

\bd{dcategoryF} {\em Define $\mathcal{F}$ to be the category of
restricted $\H_{I}$-modules $W$ of level $1$ such that $I_{W}^{\bf
p}=I$, i.e., ${\bf p}_{i}$ acts locally finitely on $W$ for every
$i\in I$.} \ed

Let $x_{i}$ $(i\in I)$ be mutually commuting independent formal
variables. Denote by $F^{0}(I,\C)$ the set of functions ${\bf \mu}:
I\rightarrow \C$ such that $\mu_{i}=0$ for all but finitely many
$i\in I$. For ${\bf \lambda}\in F^{0}(I,\C)$, set
$${\bf \lambda} {\bf x}
=\sum_{i\in I}\lambda_{i}x_{i}\in \C[x_{i}\;|\; i\in I].$$
Furthermore, set
$$M(1,{\bf \lambda)}=e^{\bf \lambda x}\C[x_{i}\;|\; i\in I],$$
a space of functions in $x_{i}$ $(i\in I)$. It is clear that
$M(1,{\bf \lambda})$ is an $\H_{I}$-module with ${\bf p}_{i}$ acting
as $\partial /\partial x_{i}$, ${\bf q}_{i}$ as (the left
multiplication of) $x_{i}$, and ${\bf k}$ as identity.

\bl{lheisenberg} The $\H_{I}$-module $M(1,{\bf \lambda})$ belongs to
the category $\mathcal{F}$ and is irreducible. For ${\bf
\lambda},{\bf \mu}\in F^{0}(I,\C)$, $M(1,{\bf \lambda})\simeq
M(1,{\bf \mu})$ if and only if ${\bf \lambda}={\bf \mu}$. \el

\begin{proof} Let $S$ be a finite subset of $I$ such that
$\lambda_{i}=0$ for all $i\in I-S$. Then $\partial e^{\bf \lambda
x}/\partial x_{i}=0$ for all $i\in I-S$. On the other hand, for any
polynomial $f({\bf x})$ in $x_{j}$ $(j\in I)$, $\partial f/\partial
x_{i}=0$ for all but finitely many $i\in I$.  It follows that for
any $w\in M(1,\lambda)$, ${\bf p}_{i}w=0$ for all but finitely many
$i\in I$. This shows that $M(1,{\bf \lambda})$ is a restricted
$\H_{I}$-module. Let $i\in I$ and let $f$ be any polynomial in
$x_{j}$ with $j\ne i$. For any $r\in \N$, $\sum_{n=0}^{r}fe^{\bf
\lambda x}\C x_{i}^{n}$ is closed under the action of ${\bf p}_{i}$.
It follows that ${\bf p}_{i}$ acts locally finitely on $M(1,{\bf
\lambda})$. This proves that $M(1,{\bf \lambda})$ belongs to the
category $\mathcal{F}$. For any function ${\bf \mu}\in F^{0}(I,\C)$,
set
$$M(1,{\bf \lambda})_{\bf \mu}=\{ w\in M(1,{\bf \lambda})\;|\; {\bf
p}_{i}w=\mu_{i}w\ \ \ \mbox{ for }i\in I\}.$$ With ${\bf p}_{i}$
acting as $\partial/\partial x_{i}$ $(i\in I)$, for $f({\bf x})\in
\C[x_{i}\;|\; i\in I]$, $e^{\bf \lambda x}f({\bf x})\in M(1,{\bf
\lambda})_{\bf \mu}$ if and only if
$$\partial f/\partial x_{i}=(\mu_{i}-\lambda_{i})f\ \ \ \mbox{ for all
}i\in I.$$ Furthermore, if $f({\bf x})\ne 0$, we have $f\in \C$ and
$\lambda_{i}=\mu_{i}$ for all $i\in I$. Thus $M(1,{\bf
\lambda})_{\bf \mu}=0$ for ${\bf \mu}\ne {\bf \lambda}$ and
$M(1,{\bf \lambda})_{\bf \lambda}=\C e^{\bf \lambda x}$. It now
follows immediately that $M(1,{\bf \lambda})$ is an irreducible
$\H_{I}$-module. {}From this proof, it is evident that for ${\bf
\lambda},{\bf \mu}\in F^{0}(I,\C)$, $M(1,{\bf \lambda})\simeq
M(1,{\bf \mu})$ if and only if ${\bf \lambda}={\bf \mu}$.
\end{proof}

\bl{levery} For every nonzero $\H_{I}$-module $W$ in the category
${\mathcal{F}}$, there exist a nonzero vector $w_{0}\in W$ and a
function $\lambda\in F^{0}(I,\C)$ such that ${\bf
p}_{i}w_{0}=\lambda_{i}w_{0}$ for $i\in I$. Furthermore, the
submodule generated by $w_{0}$ is isomorphic to $M(1,\lambda)$. \el

\begin{proof}  Let $0\ne w\in W$. {}From definition, ${\bf
p}_{i}w=0$ for $i\in I-S$, where $S$ is a finite subset of $I$. Set
$\Omega'=\{ u\in W\;|\; {\bf p}_{i}u=0\ \ \mbox{ for }i\in I-S\}$.
Then $w\in \Omega'\ne 0$. With $[{\bf p}_{r},{\bf p}_{s}]=0$ for
$r,s\in I$, $\Omega'$ is closed under the actions of ${\bf p}_{j}$
for $j\in S$. With ${\bf p}_{j}$ $(j\in S)$ locally finite and
mutually commuting, there exists $0\ne w_{0}\in \Omega'$ such that
${\bf p}_{j}w_{0}=\lambda_{j}w_{0}$ for $j\in S$ with
$\lambda_{j}\in \C$. We also have ${\bf p}_{i}w_{0}=0$ for $i\in
I-S$. Defining $\lambda_{j}=0$ for $j\in I-S$ gives rise to a
function ${\bf \lambda}\in F^{0}(I,\C)$.  In view of the P-B-W
theorem, we have $U(\H_{I})w_{0}=\C[{\bf q}_{i}\;|\; i\in I]w_{0}$.
Define a linear map
$$\psi: M(1,{\bf \lambda})\rightarrow U(\H_{I})w_{0},\ \ e^{\bf \lambda x}f({\bf
x})\mapsto f({\bf q})w_{0}.$$ It is straightforward to show that
$\psi$ is an $\H_{I}$-module isomorphism.
\end{proof}

The following is a generalization of a theorem of \cite{flm}
(Theorem 1.7.3):

\bt{theisenberg} Every irreducible $\H_{I}$-module in the category
$\mathcal{F}$ is isomorphic to $M(1,{\bf \lambda})$ for some ${\bf
\lambda}\in F^{0}(I,\C)$ and every $\H_{I}$-module in the category
$\mathcal{F}$ is completely reducible. \et

\begin{proof} The first assertion follows immediately from Lemma \ref{levery}.
For complete reducibility we first consider a special case. Let
${\mathcal{N}}$ be the subcategory of ${\mathcal{F}}$, consisting of
restricted $\H_{I}$-modules of level $1$ on which ${\bf p}_{i}$ acts
locally nilpotently for every $i\in I$. We see that for $\lambda\in
F^{0}(I,\C)$, $M(1,\lambda)$ is in ${\mathcal{N}}$ if and only if
$\lambda=0$. Then $M(1,0)$ is the only irreducible module in
${\mathcal{N}}$ up to equivalence. With Lemma \ref{levery}, the same
proof of Theorem 1.7.3 of \cite{flm} shows that every
$\H_{I}$-module in ${\mathcal{N}}$ is completely reducible.

We now prove that every $\H_{I}$-module $W$ in $\mathcal{F}$ is
completely reducible. For $\mu\in F^{0}(I,\C)$, set
$$W_{\mu}=\{ w\in W\;|\; ({\bf p}_{i}-\mu_{i})^{r}w=0
\ \ \mbox{ for }i\in I \mbox{ and for some $r\ge 0$ depending on
}i\}.$$ Note that
$$({\bf p}_{i}-\mu_{i})^{r}{\bf p}_{j}={\bf p}_{j}({\bf
p}_{i}-\mu_{i})^{r}, \ \ \ \ ({\bf p}_{i}-\mu_{i})^{r}{\bf
q}_{i}={\bf q}_{i}({\bf p}_{i}-\mu_{i})^{r}+r({\bf
p}_{i}-\mu_{i})^{r-1}{\bf k}$$ for $i,j\in I,\; r\ge 0$. It follows
that $W_{\mu}$ is an $\H_{I}$-submodule of $W$. For any $w\in W$,
since ${\bf p}_{i}w=0$ for all but finitely many $i\in I$ and since
each ${\bf p}_{i}$ is locally finite, $\C[{\bf p}_{i}\;|\; i\in I]w$
is finite-dimensional. From this we have $W=\oplus_{\mu\in
F^{0}(I,\C)}W_{\mu}$. Now it suffices to prove that for each $\mu\in
F^{0}(I,\C)$, $W_{\mu}$ is completely reducible. Let $\mu\in
F^{0}(I,\C)$ be fixed. Then ${\bf p}_{i}-\mu_{i}$ is locally
nilpotent on $W_{\mu}$ for $i\in I$. Define a linear endomorphism
$\theta_{\mu}$ of $\H_{I}$ by
$$\theta_{\mu}({\bf k})={\bf k},\ \ \theta_{\mu}({\bf p}_{i})={\bf
p}_{i}-\mu_{i} \ \mbox{ and }\ \theta_{\mu}({\bf q}_{i})={\bf q}_{i}
\ \ \ \mbox{ for }i\in I.$$ Clearly, $\theta_{\mu}$ is a Lie algebra
automorphism of $\H_{I}$. Let $\rho: \H_{I}\rightarrow \End W_{\mu}$
denote the Lie algebra homomorphism for the $\H_{I}$-module
$W_{\mu}$. Then $\rho \circ\theta_{\mu}$ is a representation of
$\H_{I}$ on $W_{\mu}$ in the category ${\mathcal{N}}$, which is
completely reducible. Consequently, $\rho$ is completely reducible.
\end{proof}

Let $I_{1}$ be any cofinite subset of $I$ and let ${\bf \lambda}\in
F^{0}(I,\C)$. Define an action of $\H_{I}$ on the space $e^{\bf
\lambda x}\C[x_{i}\;|\; i\in I]$ by
\begin{eqnarray*}
{\bf p}_{i}\mapsto \begin{cases}\partial/\partial x_{i} &
\mbox{ for }i\in I_{1}\\
x_{i} &\mbox{ for }i\in I-I_{1},
\end{cases}
\ \ \ \
{\bf q}_{i}\mapsto \begin{cases}
 x_{i} &\mbox{ for }i\in I_{1}\\
-\partial/\partial x_{i}&\mbox{ for }i\in I-I_{1}.
\end{cases}
\end{eqnarray*}
This makes $e^{\bf \lambda x}\C[x_{i}\;|\; i\in I]$ an
$\H_{I}$-module of level $1$, which we denote by
$M(1,I_{1},\lambda)$. In fact, this $\H_{I}$-module
$M(1,I_{1},\lambda)$ is a twisting of the $\H_{I}$-module $M(1,{\bf
\lambda})$ by an automorphism of $\H_{I}$. Let $\theta_{I_{1}}$ be
the linear endomorphism of $\H_{I}$ defined by
\begin{eqnarray*}
& &\theta_{I_{1}}({\bf p}_{i})={\bf p}_{i},\ \ \ \
\theta_{I_{1}}({\bf q}_{i})={\bf
q}_{i}\ \ \ \ \ \mbox{ for }i\in I_{1},\\
& &\theta_{I_{1}}({\bf p}_{i})={\bf q}_{i},\ \ \ \
\theta_{I_{1}}({\bf q}_{i})=-{\bf p}_{i}\ \ \ \mbox{ for }i\in
I-I_{1}.
\end{eqnarray*}
It is evident that $\theta_{I_{1}}$ is a Lie algebra automorphism
and that $M(1,I_{1},{\bf \lambda})$ is isomorphic to the twisting of
$M(1,{\bf \lambda})$ by the automorphism $\theta_{I_{1}}$.
Consequently, $M(1,I_{1},{\bf \lambda})$ is an irreducible
$\H_{I}$-module of level $1$. Furthermore, using the automorphism
$\theta_{I_{1}}$ and Theorem \ref{theisenberg} we immediately have:

\bc{ccategory-B} Let $W$ be an irreducible restricted
$\H_{I}$-module of level $1$ such that $I_{W}^{\bf p}\cup I_{W}^{\bf
q}=I$. Then $W$ is isomorphic to $M(1,I_{W}^{\bf p},{\bf \lambda})$
for some ${\bf \lambda}\in F^{0}(I,\C)$. \ec

Next, we continue to investigate general irreducible restricted
$\H_{I}$-modules of level $1$. Let $I=I_{0}\cup I_{1}$ be any
disjoint decomposition of $I$ with $I_{0}\ne \emptyset$ and
$I_{1}\ne \emptyset$. We view $\H_{I_{0}}$ and $\H_{I_{1}}$ as
subalgebras of $\H_{I}$ in the obvious way. Note that the two
subalgebras are commuting. Let $W_{0}$ and $W_{1}$ be irreducible
modules of level $1$ for $\H_{I_{0}}$ and $\H_{I_{1}}$,
respectively. Then $W_{0}\otimes W_{1}$ is naturally an
$\H_{I}$-module of level $1$. Furthermore, if either $I_{0}$ or
$I_{1}$ is countable, $W_{0}\otimes W_{1}$ is an irreducible
$\H_{I}$-module. (Notice that either $\H_{I_{0}}$ or $\H_{I_{1}}$ is
of countable dimension, which implies that either
$\End_{\H_{I_{0}}}(W_{0})=\C$ or $\End_{\H_{I_{1}}}(W_{1})=\C$.)

\bp{pdecomposition} Assume that $I$ is countable. Let $W$ be an
irreducible restricted $\H_{I}$-module of level $1$ such that $I\ne
I(W)=I_{W}^{\bf p}\cup I_{W}^{\bf q}$. Set $I_{0}=I-I(W)$. Then
$W\simeq M(1,I_{W}^{\bf p},{\bf \lambda})\otimes U$, where
$M(1,I_{W}^{\bf p},{\bf \lambda})$ is an $\H_{I(W)}$-module for some
${\bf \lambda}\in F^{0}(I(W),\C)$ and $U$ is an irreducible
$\H_{I_{0}}$-module such that $(I_{0})_{U}^{\bf p}=\emptyset
=(I_{0})_{U}^{\bf q}$. \ep

\begin{proof} We view $\H_{I(W)}$ and $\H_{I_{0}}$ as subalgebras of $\H_{I}$
in the obvious way. Note that the two subalgebras are commuting.
{}From Theorem \ref{theisenberg}, $W$ as an $\H_{I(W)}$-module is
completely reducible. Let $W_{1}$ be an irreducible
$\H_{I(W)}$-submodule of $W$. Since $W$ is an irreducible
$\H_{I}$-module, we have
$$W=U(\H_{I})W_{1}=U(\H_{I_{0}})W_{1}.$$
As $[\H_{I(W)},\H_{I_{0}}]=0$, it follows that $W$ as an
$\H_{I(W)}$-module is a sum of irreducible modules isomorphic to
$W_{1}$. With $I$ countable, $W_{1}$ is of countable dimension, so
that $\End_{\H_{I(W)}}(W_{1})=\C$. It follows that $W=W_{0}\otimes
W_{1}$, where $W_{0}=\Hom_{\H_{I(W)}}(W_{1},W)$ is naturally an
$\H_{I_{0}}$-module. Furthermore, $W_{0} $ is an irreducible
$\H_{I_{0}}$-module. In view of Lemma \ref{lnilp-finite} we have
$(I_{0})_{U}^{\bf p}=\emptyset =(I_{0})_{U}^{\bf q}$.
\end{proof}

Having established Proposition \ref{pdecomposition}, we next study
irreducible $\H_{I_{0}}$-modules $U$ of level $1$ with
$(I_{0})_{U}^{\bf p}=(I_{0})_{U}^{\bf q}=\emptyset$ with $I_{0}$ a
finite subset of $I$. For the rest of this section, let $I_{0}$ be a
nonempty finite subset of $I$. For ${\bf \mu}\in F^{0}(I_{0},\C)$,
set
$${\bf x^{\mu}}=\prod_{i\in I_{0}}x_{i}^{\mu_{i}},$$
where as before $x_{i}$ $(i\in I_{0})$ are mutually commuting
independent formal variables. Set
$$\C_{*}\{x_{i}\;|\; i\in I_{0}\}=\coprod_{{\bf \mu}\in F^{0}(I_{0},\C)}\C {\bf x}^{\bf \mu},$$
a vector space. With ${\bf p}_{i}$ acting as $\partial/\partial
x_{i}$ (the formal partial differential operator), ${\bf q}_{i}$ as
$x_{i}$, and ${\bf k}$ as identity, the space $\C_{*} \{ x_{i}\;|\;
i\in I_{0}\}$ becomes an $\H_{I_{0}}$-module of level $1$.

Denote by $F^{0}_{*}(I_{0},\C)$ the subset of $F^{0}(I_{0},\C)$,
consisting functions ${\bf \mu}: I_{0}\rightarrow \C$ such that
$\mu_{i}\notin \Z$ for all $i\in I_{0}$. For ${\bf \mu}\in
F^{0}_{*}(I_{0},\C)$, set
$$M_{*}[{\bf \mu}]=\coprod_{\lambda\in F(I_{0},\Z)}\C {\bf x}^{\mu+\lambda}
={\bf x^{\mu}}\C[x_{i}^{\pm 1}\;|\; i\in I_{0}],$$ where
$F(I_{0},\Z)$ denotes the set of integer-valued functions on
$I_{0}$. It is clear that $M_{*}[\bf \mu]$ is an
$\H_{I_{0}}$-submodule of $\C_{*}\{ x_{i}\;|\; i\in I_{0}\}$.

We say that an $\H_{I_{0}}$-module $W$ of level $1$ satisfies {\em
Condition $C_{I_{0}}$} if  for every $i\in I_{0}$, ${\bf q}_{i}{\bf
p}_{i}$ is semisimple and  ${\bf p}_{i}w\ne 0$, ${\bf q}_{i}w\ne 0$
for any $0\ne w\in W$. In terms of this notion we have:

\bp{pcategory-E} The $\H_{I_{0}}$-module $M_{*}[{\bf \mu}]$ is
irreducible and satisfies Condition $C_{I_{0}}$. On the other hand,
every irreducible $\H_{I_{0}}$-module $W$ satisfying Condition
$C_{I_{0}}$ is isomorphic to $M_{*}[{\bf \mu}]$ for some ${\bf
\mu}\in F^{0}_{*}(I_{0},\C))$. \ep

\begin{proof}  For any function ${\bf \beta}\in F^{0}(I_{0},\C)$,
we have $$\left(x_{i}{\partial\over\partial x_{i}}\right){\bf
x^{\beta}}=\beta_{i} {\bf x^{\beta}}\ \ \ \mbox{ for }i\in I_{0}.$$
Then ${\bf q}_{i}{\bf p}_{i}$ for all $i\in I_{0}$ act semisimply on
$M_{*}[\mu]$.  We also have
$$\frac{1}{r!}\left({\partial\over \partial x_{i}}\right)^{r}{\bf x}^{\bf \beta}
=\binom{\beta_{i}}{r}{\bf x}^{\bf \beta}x_{i}^{-r}\ \ \ \mbox{ for
}r\in \N,$$ where if $\beta_{i}\notin \Z$, $\binom{\beta_{i}}{r}\ne
0$ for $r\in \N$. It is then clear that $M_{*}[{\bf \mu}]$ is an
irreducible $\H_{I_{0}}$-module. It is also clear that
$${\bf p}_{i}w\ne 0,\ \ \ {\bf q}_{i}w\ne 0
\ \ \ \mbox{ for }i\in I_{0},\; 0\ne w\in M_{*}[\mu].$$ This proves
that $M_{*}[\mu]$ satisfies Condition $C_{I_{0}}$.

Let $W$ be an irreducible $\H_{I_{0}}$-module of level $1$,
satisfying Condition $C_{I_{0}}$. As ${\bf q}_{j}{\bf p}_{j}$ for
$j\in I_{0}$ are mutually commuting and are semisimple on $W$ by
assumption, there exists $0\ne w_{0}\in W$ such that
$${\bf q}_{j}{\bf p}_{j}w_{0}=\beta_{j}w_{0}\ \ \ \mbox{ for }j\in
I_{0},$$ where $\beta_{j}\in \C$. This gives rise to a function
${\bf \beta}\in F^{0}(I_{0},\C)$.
 We claim that ${\bf \beta}\in F^{0}_{*}(I_{0},\C)$, i.e., $\beta_{j}\notin
\Z$ for all $j\in I_{0}$.  For $j\in I_{0},\; n\in \N$, we have
\begin{eqnarray*}
& &({\bf q}_{j}{\bf p}_{j}){\bf p}_{j}^{n}w_{0}={\bf p}_{j}^{n}({\bf
q}_{j}{\bf p}_{j})w_{0}-n{\bf p}_{j}^{n}w_{0}=(\beta_{j}-n){\bf
p}_{j}^{n}w_{0},\\
& &({\bf q}_{j}{\bf p}_{j}){\bf q}_{j}^{n}w_{0}={\bf q}_{j}^{n}({\bf
q}_{j}{\bf p}_{j})w_{0}+n{\bf q}_{j}^{n}w_{0}=(\beta_{j}+n){\bf
q}_{j}^{n}w_{0}.
\end{eqnarray*}
Since ${\bf p}_{j}w\ne 0$ and ${\bf q}_{j}w\ne 0$ for any $0\ne w\in
W$, it follows that $\beta_{j}-n\ne 0$ and $\beta_{j}+n\ne 0$ for
$n\ge 0$, proving $\beta_{j}\notin \Z$. Thus, $\beta\in
F_{*}^{0}(I_{0},\C)$.

Since $W$ is an irreducible $\H_{I_{0}}$-module, we have
$W=U(\H_{I_{0}})w_{0}$. Define a linear map $\psi: M_{*}[{\bf
\beta}]\rightarrow W$ by
$$\psi({\bf x}^{\bf \beta+m})={\bf q}^{\bf m}w_{0}$$
for ${\bf m}: I_{0}\rightarrow \Z$, where ${\bf q}^{\bf
m}=\prod_{i\in I_{0}}{\bf q}_{i}^{m_{i}}$ with ${\bf
q}_{i}^{m_{i}}={\bf q}_{i}^{m_{i}}$ for $m_{i}\ge 0$ and
$${\bf
q}_{i}^{m_{i}}=\frac{1}{\binom{\beta_{i}}{-m_{i}}}\frac{1}{(-m_{i})!}{\bf
p}_{i}^{-m_{i}}$$ for $m_{i}<0$. One can show that $\psi$ is an
$\H_{I_{0}}$-module isomorphism.
\end{proof}

\section{Vertex algebras $V_{(\frak{h},L)}$ and their modules}
In this section, we study vertex algebras associated with possibly
degenerate even lattices. This slightly generalizes the vertex
algebras associated with nondegenerate even lattices. We construct
and classify irreducible modules satisfying certain conditions for
the vertex algebras.

First, we start with vertex operator algebras associated with
(infinite-dimensional) Heisenberg Lie algebras. Let $\frak{h}$ be a
finite-dimensional vector space equipped with a nondegenerate
symmetric bilinear form $\<\cdot,\cdot\>$, which is fixed throughout
this section. Viewing $\frak{h}$ as an abelian Lie algebra equipped
with $\<\cdot,\cdot\>$ as a nondegenerate symmetric invariant
bilinear form, we have an affine Lie algebra
$$\hat{\frak{h}}=\frak{h}\otimes \C[t,t^{-1}]\oplus \C {\bf k},$$
where ${\bf k}$ is central and for $u,v\in \frak{h},\; m,n\in \Z$,
\begin{eqnarray}
[u(m),v(n)]=m\<u,v\>\delta_{m+n,0}{\bf k},
\end{eqnarray}
where $h(n)=h\otimes t^{n}$ for $h\in \frak{h},\; n\in \Z$. For
$n\in \Z$, set
$$\frak{h}(n)=\{ h(n)\;|\; h\in \frak{h}\}\subset \hat{\frak{h}}$$
 and we set
 $$\hat{\frak{h}}_{*}=\coprod_{n\ne 0}\frak{h}(n) +\C {\bf k}.$$
Note that $\frak{h}(0)$ is a central subalgebra,
$\hat{\frak{h}}_{*}$ is a Heisenberg algebra, and
$$\hat{\frak{h}}=\hat{\frak{h}}_{*}\oplus \frak{h}(0),$$
a Lie algebra product decomposition.

\bl{lnilp} Let $W$ be any irreducible $\hat{\frak{h}}$-module. Then
 ${\bf k}$ and $h(0)$ for every $h\in \frak{h}$ act as scalars on
$W$ and $W$ as an $\hat{\frak{h}}_{*}$-module is also irreducible.
 \el

\begin{proof} Since $W$ is an irreducible $\hat{\frak{h}}$-module and
$\frak{h}$ is finite-dimensional, $W$ is of countable dimension over
$\C$. By a version of Schur lemma, we have $\End_{\hat{\frak{h}}}
(W)=\C$. With ${\bf k}$ and $h(0)$ being central in
$\hat{\frak{h}}$, they must act as scalars. It is now clear that $W$
as an $\hat{\frak{h}}_{*}$-module is also irreducible.
\end{proof}

An $\hat{\frak{h}}$-module $W$ is said to be of {\em level} $\ell$
in $\C$ if ${\bf k}$ acts on $W$ as scalar $\ell$, and an
$\hat{\frak{h}}$-module $W$ is said to be {\em restricted} if for
every $w\in W$, $\frak{h}(n)w=0$ for $n$ sufficiently large. For any
$\ell\in \C$, denote by $\C_{\ell}$ the $1$-dimensional
$\frak{h}\otimes \C[t]+\C {\bf k}$-module $\C$ with $\frak{h}\otimes
\C[t]$ acting trivially and with ${\bf k}$ acting as scalar $\ell$.
Form the induced module
$$V_{\hat{\frak{h}}}(\ell,0)
=U(\hat{\frak{h}})\otimes_{U(\frak{h}\otimes \C[t]+\C {\bf k})}
\C_{\ell}.$$ Set ${\bf 1}=1\otimes 1\in V_{\hat{\frak{h}}}(\ell,0)$
and identify $\frak{h}$ as a subspace of
$V_{\hat{\frak{h}}}(\ell,0)$ through the linear map $h\mapsto
h(-1){\bf 1}$. It is well known now (cf. \cite{ll}) that there
exists a unique vertex algebra structure on
$V_{\hat{\frak{h}}}(\ell,0)$ with ${\bf 1}$ as the vacuum vector and
with $Y(h,x)=h(x)=\sum_{n\in \Z}h(n)x^{-n-1}$ for $h\in \frak{h}$.
Furthermore, for every nonzero $\ell$, $V_{\hat{\frak{h}}}(\ell,0)$
is a vertex operator algebra of central charge $d=\dim \frak{h}$
with conformal vector
$$\omega=\frac{1}{2\ell}\sum_{r=1}^{d}h^{(r)}(-1)h^{(r)}(-1){\bf 1},$$
where $\{h^{(1)},\dots,h^{(d)}\}$ is any orthonormal basis of
$\frak{h}$.  It is also known (cf. \cite{ll}) that every module
$(W,Y_{W})$ for vertex algebra $V_{\hat{\frak{h}}}(\ell,0)$ is a
restricted $\hat{\frak{h}}$-module of level $\ell$ with $h(x)$
acting as $Y_{W}(h,x)$ for $h\in \frak{h}$ and the set of
$V_{\hat{\frak{h}}}(\ell,0)$-submodules of $W$ coincides with the
set of $\hat{\frak{h}}$-submodules of $W$.  On the other hand, on
every restricted $\hat{\frak{h}}$-module $W$ of level $\ell$, there
is a unique module structure $Y_{W}$ for vertex algebra
$V_{\hat{\frak{h}}}(\ell,0)$ with $Y_{W}(h,x)=h(x)$ for $h\in
\frak{h}$.

It is known (cf. \cite{ll}) that vertex algebras
$V_{\hat{\frak{h}}}(\ell,0)$ for $\ell\ne 0$ are all isomorphic. In
view of this, we restrict ourselves to the vertex operator algebra
$V_{\hat{\frak{h}}}(1,0)$. In literature, the vertex operator
algebra $V_{\hat{\frak{h}}}(1,0)$ is also often denoted by $M(1)$.
In view of Lemma \ref{lnilp}, classifying irreducible modules for
vertex algebra $M(1)$ amounts to classifying irreducible restricted
$\hat{\frak{h}}_{*}$-modules of level $1$. One way to apply the
results of Section 2 is to fix an orthonormal basis
$\{h^{(1)},\dots,h^{(d)}\}$ of $\frak{h}$ and set
$$I=\{(r,n)\;|\; 1\le r\le d,\; n\ge 1\}.$$
Then identify the Heisenberg Lie algebra $\hat{\frak{h}}_{*}$ with
$\H_{I}$ by
$${\bf p}_{(r,n)}=\frac{1}{n}h^{(r)}(n),\ \ {\bf
q}_{(r,n)}=h^{(r)}(-n)\ \ \ \mbox{ for }(r,n)\in I.$$

In view of Theorem \ref{theisenberg} we immediately have:

\bp{pcomplete} Let $W$ be a restricted $\hat{\frak{h}}$-module of
level $1$ such that $\frak{h}(0)$ is semisimple and $\frak{h}(n)$ is
locally finite for $n\ge 1$. Then $W$ is completely reducible. \ep

For $\alpha\in \frak{h}$, we set (\cite{lw}, \cite{flm})
\begin{eqnarray}
E^{\pm}(\alpha,x) =\exp \left(\sum_{n=1}^{\infty}\frac{\alpha (\pm
n)}{\pm n}x^{\mp n}\right).
\end{eqnarray}
The following are the fundamental properties:
\begin{eqnarray*}
& &E^{\pm}(0,x)=1,\\
& &E^{\pm}(\alpha,x_{1})E^{\pm }(\beta,x_{2})= E^{\pm
}(\beta,x_{2})E^{\pm}(\alpha,x_{1}),\\
& &E^{+}(\alpha,x_{1})E^{-}(\beta,x_{2})
=\left(1-\frac{x_{2}}{x_{1}}\right)^{\<\alpha,\beta\>}E^{-}(\beta,x_{2})E^{+}(\alpha,x_{1}),\\
& &E^{\pm}(\alpha,x)E^{\pm }(\beta,x)= E^{\pm }(\alpha+\beta,x).
\end{eqnarray*}

Set
\begin{eqnarray*}
\bar{\Delta}(\alpha,x)=(-x)^{\alpha(0)}\exp\left(\sum_{n=1}^{\infty}\frac{\alpha(n)}{-n}
(-x)^{-n}\right)=(-x)^{\alpha(0)}E^{+}(-\alpha,-x).
\end{eqnarray*}
This is a well defined element of $(\End W)[[x,x^{-1}]]$ for any
module $W$ for vertex algebra $M(1)$, on which $\alpha(0)$ acts
semisimply with only integer eigenvalues and $\alpha(n)$ for $n\ge
1$ act locally nilpotently. The following are immediate
consequences:
\begin{eqnarray*}
& &\bar{\Delta}(0,x)=1,\\
 & &\bar{\Delta}(\alpha,x_{1})\bar{\Delta}(\beta,x_{2})
=\bar{\Delta}(\beta,x_{2})\bar{\Delta}(\alpha,x_{1}),\\
& &\bar{\Delta}(\alpha,x)\bar{\Delta}(\beta,x)
=\bar{\Delta}(\alpha+\beta,x).
\end{eqnarray*}

\br{rrecall-property} {\em Let $\alpha\in \frak{h}$ and let
$(W,Y_{W})$ be an $M(1)$-module on which  $\alpha(0)$ acts
semisimply with only integer eigenvalues and $\alpha(n)$ for $n\ge
1$ act locally nilpotently. Recall from \cite{li-phys} (and
\cite{li-ext}) that
\begin{eqnarray*}
\Delta(\alpha,x)=x^{\alpha(0)}E^{+}(-\alpha,-x).
\end{eqnarray*}
The following basic properties were established in \cite{li-phys}
and \cite{li-ext}\footnote{Note that the notation $E^{-}(\alpha,x)$
in \cite{li-ext} is the notation $E^{-}(-\alpha,x)$ in this paper
and in \cite{flm}, and  all the related formulas have been
 adjusted correspondingly.}: $\Delta(\alpha,x){\bf 1}={\bf 1}$,
\begin{eqnarray}
\Delta(\alpha,x_{1})Y_{W}(v,x_{2})=Y_{W}(\Delta(\alpha,x_{1}+x_{2})v,x_{2})\Delta(\alpha,x_{1})
\end{eqnarray}
for $v\in M(1)$, and
\begin{eqnarray*}
& &Y_{W}(E^{-}(\alpha,x_{1})v,x_{2})\\
&=&E^{-}(\alpha,x_{1}+x_{2})E^{-}(-\alpha,x_{2})Y_{W}(v,x_{2})
x_{2}^{\alpha(0)}E^{+}(-\alpha,x_{2})(x_{2}+x_{1})^{-\alpha(0)}E^{+}(\alpha,x_{2}+x_{1})\\
&=&E^{-}(\alpha,x_{1}+x_{2})E^{-}(-\alpha,x_{2})Y_{W}(v,x_{2})
\Delta(\alpha,-x_{2})\Delta(-\alpha,-x_{2}-x_{1}),
\end{eqnarray*}
\begin{eqnarray}\label{e2.4}
E^{-}(\alpha,x_{1})Y_{W}(v,x_{2})E^{-}(-\alpha,x_{1})
=Y_{W}(\Delta(\alpha,x_{2}-x_{1})\Delta(-\alpha,x_{2})v,x_{2})
\end{eqnarray}
for $\alpha\in \frak{h},\; v\in M(1)$.  Note that
$\bar{\Delta}(\alpha,x)=(-1)^{\alpha(0)}\Delta(\alpha,x)$ and that
$$(-1)^{\alpha(0)}Y_{W}(v,x)=Y_{W}((-1)^{\alpha(0)}v,x)(-1)^{\alpha(0)}
\ \ \ \mbox{ for }v\in M(1).$$ It follows that all the above
properties with $\Delta(\alpha,x)$ replaced by
$\Bar{\Delta}(\alpha,x)$ still hold.} \er

For the rest of this section, we assume that $L$ is an additive
subgroup of $\frak{h}$ such that
\begin{eqnarray}\label{eevenL}
\<\alpha,\alpha\>\in 2\Z \ \ \ \mbox{ for }\alpha \in L,
\end{eqnarray}
equipped with a function $\varepsilon: L\times L\rightarrow
\C^{\times}$, satisfying the condition that
\begin{eqnarray}
& &\varepsilon(\alpha,0)=\varepsilon(0,\alpha)=1,\label{enormalization}\\
& &\varepsilon(\alpha,\beta+\gamma)\varepsilon(\beta,\gamma)
=\varepsilon(\alpha+\beta,\gamma)\varepsilon(\alpha,\beta),\\
& &\varepsilon(\alpha,\beta)\varepsilon(\beta,\alpha)^{-1}
=(-1)^{\<\alpha,\beta\>}\label{e2cocycle}
\end{eqnarray}
for $\alpha,\beta,\gamma\in L$.

 \br{rcocycle-exist} {\em
Note that if $L$ is finitely generated, $L$ is a finite rank free
group as $L$ is torsion-free. Assume that $L$ is a free group of
finite rank. Let $\{\alpha_{1},\dots,\alpha_{r}\}$ be a $\Z$-basis
of $L$. Define $\varepsilon: L\times L\rightarrow \C^{\times}$ to be
the group homomorphism uniquely determined by
\begin{eqnarray*}
\varepsilon(\alpha_{i},\alpha_{j})
=\begin{cases}(-1)^{\<\alpha_{i},\alpha_{j}\>}\ \ &\mbox{ if }i<j\\
1 & \mbox{ if }i\ge j.
\end{cases}
\end{eqnarray*}
Then $\varepsilon$ satisfies all the conditions. } \er

Set
\begin{eqnarray}
V_{(\frak{h},L)}=\C[L]\otimes M(1).
\end{eqnarray}
For $\alpha,\beta\in L,\; u,v\in M(1)$, we define
\begin{eqnarray}\label{eexplicitformula1}
& &Y(e^{\alpha}\otimes u,x)(e^{\beta}\otimes v)\nonumber\\
&=&\varepsilon(\alpha,\beta)e^{\alpha+\beta}\otimes
x^{\<\alpha,\beta\>}E^{-}(-\alpha,x)Y(\bar{\Delta}(\beta,x)u,x)\bar{\Delta}(\alpha,-x)v.
\end{eqnarray}
In particular, for $h\in \frak{h}$,
\begin{eqnarray}
& &Y(e^{0}\otimes h,x)(e^{\beta}\otimes v)=e^{\beta}\otimes
(\<\beta,h\>x^{-1}+Y(h,x))v,
\label{eY-formula-special-h}\\
 & & Y(e^{\alpha}\otimes {\bf
1},x)(e^{\beta}\otimes
v)=\varepsilon(\alpha,\beta)e^{\alpha+\beta}\otimes
x^{\<\alpha,\beta\>}E^{-}(-\alpha,x)E^{+}(-\alpha,x)v\ \ \ \ \
\label{eY-formula-special}
\end{eqnarray}
as $\bar{\Delta}(\beta,x)h=h+\<\beta,h\>{\bf 1}x^{-1}$ and
$\alpha(0)=0$ on $M(1)$. One can prove that the quadruple
$(V_{(\frak{h},L)},Y, e^{0}\otimes {\bf 1})$ carries the structure
of a vertex algebra, by using a theorem of \cite{fkrw} and
\cite{mp}. Here, we give a uniform treatment for both the vertex
algebras and their modules.

\bd{dconditionCL} {\em We say an $\hat{\frak{h}}$-module $W$
satisfies {\em Condition $C_{L}$} if $W$ is a restricted
$\hat{\frak{h}}$-module of level $1$, satisfying the condition that
for any $\alpha\in L$, $\alpha(0)$ is semisimple with only integer
eigenvalues and $\alpha(n)$ with $n\ge 1$ are locally nilpotent.}
\ed

\bt{tmain} Let $U$ be an $\hat{\frak{h}}$-module of level $1$,
satisfying Condition $C_{L}$. Set
\begin{eqnarray}
V_{(\frak{h},L)}(U)=\C[L]\otimes U.
\end{eqnarray}
For $\alpha,\beta\in L,\; v\in M(1),\; w\in U$,
we define
\begin{eqnarray}\label{eformula-def}
& &Y_{W}(e^{\alpha}\otimes v,x)(e^{\beta}\otimes
w)\nonumber\\
&=&\varepsilon(\alpha,\beta)e^{\alpha+\beta}\otimes
x^{\<\alpha,\beta\>}E^{-}(-\alpha,x)Y_{U}(\bar{\Delta}(\beta,x)v,x)\bar{\Delta}(\alpha,-x)w.
\end{eqnarray}
Then $(V_{(\frak{h},L)},Y, e^{0}\otimes {\bf 1})$ carries the
structure of a vertex algebra with $M(1)$ as a vertex subalgebra and
$(V_{(\frak{h},L)}(U),Y_{W})$ carries the structure of a
$V_{(\frak{h},L)}$-module. Furthermore, if $U$ is irreducible, so is
$V_{(\frak{h},L)}(U)$. \et

\begin{proof} First, notice that the $\hat{\frak{h}}$-module $M(1)$ satisfies
all the assumptions on $U$ and that
$V_{(\frak{h},L)}=V_{(\frak{h},L)}(U)$ with $U=M(1)$ where $Y_{W}$
coincides with $Y$. Second, for any $\alpha,\beta\in L,\; v\in
M(1),\; w\in U$, we have
$$Y_{W}(e^{\alpha}\otimes v,x)(e^{\beta}\otimes w)\in V_{(\frak{h},L)}(U)((x)).$$
Third, $Y_{W}(e^{0}\otimes {\bf 1},x)=1$ and when $U=M(1)$ we also
have
$$Y(e^{\alpha}\otimes u,x)(e^{0}\otimes {\bf 1})
=e^{\alpha}\otimes E^{-}(-\alpha,x)Y(u,x){\bf 1}\in
V_{(\frak{h},L)}[[x]]$$ with $$\lim_{x\rightarrow
0}Y(e^{\alpha}\otimes u,x)(e^{0}\otimes {\bf 1})=e^{\alpha}\otimes
u.$$
 Next, we show
that the Jacobi identity holds. Let $\alpha,\beta,\gamma\in L,\;
u,v\in M(1),\; w\in U$. Using definition and formulas we have
\begin{eqnarray*}
& &Y_{W}(e^{\alpha}\otimes u,x_{1})Y_{W}(e^{\beta}\otimes
v,x_{2})(e^{\gamma}\otimes w)\\
&=&Y_{W}(e^{\alpha}\otimes
u,x_{1})\left(\varepsilon(\beta,\gamma)e^{\beta+\gamma}\otimes
x_{2}^{\<\beta,\gamma\>}
E^{-}(-\beta,x_{2})Y_{U}(\bar{\Delta}(\gamma,x_{2})v,x_{2})\bar{\Delta}(\beta,-x_{2})w\right)\\
 &=&\varepsilon(\alpha,\beta+\gamma)\varepsilon(\beta,\gamma)e^{\alpha+\beta+\gamma}\otimes
 x_{1}^{\<\alpha,\beta+\gamma\>}x_{2}^{\<\beta,\gamma\>}
 E^{-}(-\alpha,x_{1})\cdot\\
 & &\cdot Y_{U}(\bar{\Delta}(\beta+\gamma,x_{1})u,x_{1})\bar{\Delta}(\alpha,-x_{1})
 E^{-}(-\beta,x_{2})Y_{U}(\bar{\Delta}(\gamma,x_{2})v,x_{2})\bar{\Delta}(\beta,-x_{2})w\\
&=&\varepsilon(\alpha,\beta+\gamma)\varepsilon(\beta,\gamma)e^{\alpha+\beta+\gamma}\otimes
 x_{1}^{\<\alpha,\beta+\gamma\>}x_{2}^{\<\beta,\gamma\>}
 E^{-}(-\alpha,x_{1}) Y_{U}(\bar{\Delta}(\beta+\gamma,x_{1})u,x_{1})\cdot\\
 & &\cdot (1-x_{2}/x_{1})^{\<\alpha,\beta\>}
 E^{-}(-\beta,x_{2})Y_{U}(\bar{\Delta}(\alpha,-x_{1}+x_{2})\bar{\Delta}(\gamma,x_{2})v,x_{2})
 \bar{\Delta}(\alpha,-x_{1})\bar{\Delta}(\beta,-x_{2})w\\
\end{eqnarray*}
By (\ref{e2.4}) we have
\begin{eqnarray*}
& &Y_{U}(\bar{\Delta}(\beta+\gamma,x_{1})u,x_{1})E^{-}(-\beta,x_{2})\\
&=&E^{-}(-\beta,x_{2})
Y_{U}(\bar{\Delta}(\beta,x_{1}-x_{2})\bar{\Delta}(-\beta,x_{1})\bar{\Delta}(\beta+\gamma,x_{1})u,x_{1})\\
&=&E^{-}(-\beta,x_{2})
Y_{U}(\bar{\Delta}(\beta,x_{1}-x_{2})\bar{\Delta}(\gamma,x_{1})u,x_{1}).
\end{eqnarray*}
Then
\begin{eqnarray*}
& &Y_{W}(e^{\alpha}\otimes u,x_{1})Y_{W}(e^{\beta}\otimes
v,x_{2})(e^{\gamma}\otimes w)\\
&=&\varepsilon(\alpha,\beta+\gamma)\varepsilon(\beta,\gamma)e^{\alpha+\beta+\gamma}\otimes
 x_{1}^{\<\alpha,\gamma\>}x_{2}^{\<\beta,\gamma\>}(x_{1}-x_{2})^{\<\alpha,\beta\>}
 E^{-}(-\alpha,x_{1})E^{-}(-\beta,x_{2})\cdot\\
& &\cdot
Y_{U}(\bar{\Delta}(\beta,x_{1}-x_{2})\bar{\Delta}(\gamma,x_{1})u,x_{1})
Y_{U}(\bar{\Delta}(\alpha,-x_{1}+x_{2})\bar{\Delta}(\gamma,x_{2})v,x_{2})
 \bar{\Delta}(\alpha,-x_{1})\bar{\Delta}(\beta,-x_{2})w.
\end{eqnarray*}
This also shows
\begin{eqnarray*}
& &Y_{W}(e^{\beta}\otimes v,x_{2})Y_{W}(e^{\alpha}\otimes u,x_{1})
(e^{\gamma}\otimes w)\\
&=&\varepsilon(\beta,\alpha+\gamma)\varepsilon(\alpha,\gamma)e^{\alpha+\beta+\gamma}\otimes
 x_{1}^{\<\alpha,\gamma\>}x_{2}^{\<\beta,\gamma\>}(x_{2}-x_{1})^{\<\beta,\alpha\>}
 E^{-}(-\alpha,x_{1})E^{-}(-\beta,x_{2})\cdot\\
& &\cdot
Y_{U}(\bar{\Delta}(\alpha,x_{2}-x_{1})\bar{\Delta}(\gamma,x_{2})v,x_{2})
Y_{U}(\bar{\Delta}(\beta,-x_{2}+x_{1})\bar{\Delta}(\gamma,x_{1})u,x_{1})
\bar{\Delta}(\alpha,-x_{1})\bar{\Delta}(\beta,-x_{2})w\\
&=&\varepsilon(\alpha,\beta)\varepsilon(\alpha+\beta,\gamma)e^{\alpha+\beta+\gamma}\otimes
 x_{1}^{\<\alpha,\gamma\>}x_{2}^{\<\beta,\gamma\>}(-x_{2}+x_{1})^{\<\beta,\alpha\>}
 E^{-}(-\alpha,x_{1})E^{-}(-\beta,x_{2})\cdot\\
& &\cdot
Y_{U}(\bar{\Delta}(\alpha,x_{2}-x_{1})\bar{\Delta}(\gamma,x_{2})v,x_{2})
Y_{U}(\bar{\Delta}(\beta,-x_{2}+x_{1})\bar{\Delta}(\gamma,x_{1})u,x_{1})
 \bar{\Delta}(\alpha,-x_{1})\bar{\Delta}(\beta,-x_{2})w,
\end{eqnarray*}
where $\varepsilon(\beta,\alpha+\gamma)\varepsilon(\alpha,\gamma)
=(-1)^{\<\alpha,\beta\>}\varepsilon(\alpha,\beta)\varepsilon(\alpha+\beta,\gamma)$.
On the other hand, we have
\begin{eqnarray*}
& &Y_{W}\left(Y(e^{\alpha}\otimes u,x_{0})(e^{\beta}\otimes
v),x_{2}\right)(e^{\gamma}\otimes w)\\
&=&Y_{W}(\varepsilon(\alpha,\beta)e^{\alpha+\beta}\otimes
x_{0}^{\<\alpha,\beta\>}
E^{-}(-\alpha,x_{0})Y(\bar{\Delta}(\beta,x_{0})u,x_{0})\bar{\Delta}(\alpha,-x_{0})v,x_{2})
(e^{\gamma}\otimes w)\\
&=&\varepsilon(\alpha+\beta,\gamma)\varepsilon(\alpha,\beta)e^{\alpha+\beta+\gamma}\otimes
x_{0}^{\<\alpha,\beta\>}x_{2}^{\<\alpha+\beta,\gamma\>}E^{-}(-\alpha-\beta,x_{2})\cdot\\
& &\cdot Y_{U}\left(\bar{\Delta}(\gamma,x_{2})E^{-}(-\alpha,x_{0})
Y(\bar{\Delta}(\beta,x_{0})u,x_{0})\bar{\Delta}(\alpha,-x_{0})v,x_{2}\right)
\bar{\Delta}(\alpha+\beta,-x_{2})w\\
&=&\varepsilon(\alpha+\beta,\gamma)\varepsilon(\alpha,\beta)e^{\alpha+\beta+\gamma}\otimes
x_{0}^{\<\alpha,\beta\>}x_{2}^{\<\alpha+\beta,\gamma\>}E^{-}(-\alpha-\beta,x_{2})
(1+x_{0}/x_{2})^{\<\alpha,\gamma\>}\cdot\\
& &\cdot Y_{U}\left(E^{-}(-\alpha,x_{0})\Delta(\gamma,x_{2})
Y(\bar{\Delta}(\beta,x_{0})u,x_{0})\bar{\Delta}(\alpha,-x_{0})v,x_{2}\right)
\bar{\Delta}(\alpha+\beta,-x_{2})w\\
&=&\varepsilon(\alpha+\beta,\gamma)\varepsilon(\alpha,\beta)e^{\alpha+\beta+\gamma}\otimes
x_{0}^{\<\alpha,\beta\>}x_{2}^{\<\beta,\gamma\>}E^{-}(-\alpha-\beta,x_{2})
(x_{2}+x_{0})^{\<\alpha,\gamma\>}\cdot\\
& &\cdot Y_{U}\left(E^{-}(-\alpha,x_{0})
Y(\bar{\Delta}(\gamma,x_{2}+x_{0})\bar{\Delta}(\beta,x_{0})u,x_{0})
\bar{\Delta}(\gamma,x_{2})\bar{\Delta}(\alpha,-x_{0})v,x_{2}\right)
\cdot\\
& &\ \ \ \ \cdot \bar{\Delta}(\alpha+\beta,-x_{2})w\\
&=&\varepsilon(\alpha+\beta,\gamma)\varepsilon(\alpha,\beta)e^{\alpha+\beta+\gamma}\otimes
x_{0}^{\<\alpha,\beta\>}x_{2}^{\<\beta,\gamma\>}E^{-}(-\alpha-\beta,x_{2})
(x_{2}+x_{0})^{\<\alpha,\gamma\>}
E^{-}(-\alpha,x_{0}+x_{2})\cdot\\
& &\cdot E^{-}(\alpha,x_{2})Y_{U}\left(
Y(\bar{\Delta}(\gamma,x_{2}+x_{0})\bar{\Delta}(\beta,x_{0})u,x_{0})
\bar{\Delta}(\gamma,x_{2})\bar{\Delta}(\alpha,-x_{0})v,x_{2}\right)
\cdot\\
& &\ \ \ \ \cdot
\bar{\Delta}(-\alpha,-x_{2})\bar{\Delta}(\alpha,-x_{2}-x_{0})
\bar{\Delta}(\alpha+\beta,-x_{2})w\\
&=&\varepsilon(\alpha+\beta,\gamma)\varepsilon(\alpha,\beta)e^{\alpha+\beta+\gamma}\otimes
x_{0}^{\<\alpha,\beta\>}x_{2}^{\<\beta,\gamma\>}(x_{2}+x_{0})^{\<\alpha,\gamma\>}E^{-}(-\beta,x_{2})
 E^{-}(-\alpha,x_{0}+x_{2})\cdot\\
& &\cdot Y_{U}\left(
Y(\bar{\Delta}(\gamma,x_{2}+x_{0})\bar{\Delta}(\beta,x_{0})u,x_{0})
\bar{\Delta}(\gamma,x_{2})\bar{\Delta}(\alpha,-x_{0})v,x_{2}\right)
\cdot\\
& &\ \ \ \ \cdot \bar{\Delta}(\alpha,-x_{2}-x_{0})
\bar{\Delta}(\beta,-x_{2})w.
\end{eqnarray*}
Set
\begin{eqnarray*} & &A=\bar{\Delta}(\gamma,x_{1})\bar{\Delta}(\beta,x_{0})u, \
\ B=\bar{\Delta}(\gamma,x_{2})\bar{\Delta}(\alpha,-x_{0})v\in
M(1)[x_{1}^{\pm
1},x_{2}^{\pm 1},x_{0}^{\pm 1}],\\
& &\ \ \ \ \ C=\bar{\Delta}(\alpha,-x_{1})
\bar{\Delta}(\beta,-x_{2})w\in U[x_{1}^{\pm 1},x_{2}^{\pm 1}].
\end{eqnarray*}
We have the following Jacobi identity
\begin{eqnarray*}
& &x_{0}^{-1}\delta\left(\frac{x_{1}-x_{2}}{x_{0}}\right)
Y_{U}(A,x_{1})Y_{U}(B,x_{2})C
-x_{0}^{-1}\delta\left(\frac{x_{2}-x_{1}}{-x_{0}}\right)
Y_{U}(B,x_{2})Y_{U}(A,x_{1})C\ \ \\
&&\ \ \ =x_{1}^{-1}\delta\left(\frac{x_{2}+x_{0}}{x_{1}}\right)
Y_{U}(Y(A,x_{0})B,x_{2})C.
\end{eqnarray*}
Using this and delta-function substitutions we obtain the Jacobi
identity as desired. This proves the assertions on vertex algebra
structure and module structure. It follows from
(\ref{eY-formula-special-h}) that $M(1)$ is a vertex subalgebra with
$e^{0}\otimes v$ identified with $v$ for $v\in M(1)$. Furthermore,
by (\ref{eformula-def}) (cf. (\ref{eY-formula-special-h})) we have
\begin{eqnarray}
& &h_{W}(n)(e^{\beta}\otimes w)=e^{\beta}\otimes h(n)w\ \ \ \mbox{ for }n\ne 0,\\
& &h_{W}(0)(e^{\beta}\otimes w)=e^{\beta}\otimes (\<\beta,h\>+h(0))w
\end{eqnarray}
for $h\in \frak{h},\; \alpha,\beta\in L$, where
$Y_{W}(h,x)=\sum_{n\in \Z}h_{W}(n)x^{-n-1}$.

Assume that $U$ is an irreducible $\hat{\frak{h}}$-module. In view
of Lemma \ref{lnilp}, there exists $\lambda\in \frak{h}$ such that
$h(0)$ acts as scalar $\<\lambda,h\>$ on $U$ for $h\in \frak{h}$.
Furthermore, $h(0)$ acts on $\C e^{\alpha}\otimes U$ as scalar
$\<\lambda+\alpha,h\>$ for $\alpha\in L$. Consequently, $\C
e^{\alpha}\otimes U$ for $\alpha\in L$ are non-isomorphic
irreducible $\hat{\frak{h}}$-submodules of $V_{(\frak{h},L)}(U)$.
Furthermore, for $\alpha,\beta\in L,\; w\in U$ we have
$$Y_{W}(e^{\alpha},x)(e^{\beta}\otimes w)
=\varepsilon(\alpha,\beta)e^{\alpha+\beta}\otimes
x^{\<\alpha,\beta\>}E^{-}(-\alpha,x)E^{+}(-\alpha,x)x^{\alpha(0)}w,$$
which gives
$$x^{-\alpha(0)}E^{+}(\alpha,x)\left(E^{-}(\alpha,x)Y_{W}(e^{\alpha},x)(e^{\beta}\otimes
w)\right) =\varepsilon(\alpha,\beta)e^{\alpha+\beta}\otimes
x^{\<\alpha,\beta\>}w.$$ Then for every $\beta\in L$, $0\ne w\in U$,
$e^{\beta}\otimes w$ generates $V_{(\frak{h},L)}(U)$ as an
$V_{(\frak{h},L)}$-module. It now follows that $V_{(\frak{h},L)}(U)$
is an irreducible $V_{(\frak{h},L)}$-module.
\end{proof}

\br{rlx} {\em  Let $V$ be a simple vertex algebra containing $M(1)$
(associated to some $\frak{h}$) as a vertex subalgebra  such that
$V$ as an $M(1)$-module is a direct sum of non-isomorphic
irreducible modules $M(1)\otimes \C e^{\alpha}$ with $\alpha\in S$,
where $S$ is a (nonempty) subset of $\frak{h}$. It was proved in
\cite{lx} that $S$ is an (additive) subgroup equipped with a
function $\varepsilon: S\times S\rightarrow \C^{\times}$, satisfying
all the conditions (\ref{eevenL})---(\ref{e2cocycle}), such that
(\ref{eY-formula-special}) holds for $\alpha,\beta\in S$. Now,
Theorem \ref{tmain} provides the existence of such a simple vertex
algebra. On the other hand, it is straightforward to see that
cohomologous cocycles $\varepsilon$ give rise to isomorphic vertex
algebras. If we restrict ourselves to finitely generated vertex
algebras, then $S$ is finitely generated. As $S$, a subgroup of
$\frak{h}$, is torsion-free, $S$ is a finite rank free group. It
follows from \cite{flm} (Chapter 5) that all the $2$-cocycles
$\varepsilon: S\times S\rightarrow \C^{\times}$ satisfying
(\ref{enormalization})---(\ref{e2cocycle}) are cohomologous. }\er

 \bl{lconverse} Let $W$ be any irreducible
$V_{(\frak{h},L)}$-module. Then $W$ is an $\hat{\frak{h}}$-module
with $h(x)=Y_{W}(h,x)$ for $h\in \frak{h}$, satisfying Condition
$C_{L}$. \el

\begin{proof} With $M(1)$ as a vertex subalgebra of $V_{(\frak{h},L)}$, $W$ is
a module for vertex algebra $M(1)$, so that $W$ is a restricted
$\hat{\frak{h}}$-module of level $1$. By Lemma 3.15 of
\cite{dlm-reg}, there exists a nonzero vector $w\in W$ such that
$$\alpha(n)w=0\ \ \ \ \mbox{ for }\alpha\in L,\; n\ge 1.$$
For the rest of the proof we follow Dong's arguments in
\cite{dong1}. Let $\alpha$ be any nonzero element of $L$. As
$V_{(\frak{h},L)}$ is simple, from [DL], we have
$Y(e^{\alpha},x)w\ne 0$. Assume that $e^{\alpha}_{k}w\ne 0$ and
$e^{\alpha}_{m}w=0$ for $m>k$. From the relation ${d\over dx}
e^{\alpha}(x)w=\alpha(x)e^{\alpha}(x)w$, extracting the coefficients
of $x^{-k-2}$ we obtain
$(-k-1)e^{\alpha}_{k}w=\alpha(0)e^{\alpha}_{k}w$. That is,
$e^{\alpha}_{k}w$ is an eigenvector of $\alpha(0)$ with integer
eigenvalue $-k-1$. We have
$$[\alpha(m),h(n)]=m\<\alpha,h\>\delta_{m+n,0},\ \ \
[\alpha(m),e^{\beta}_{n}]=\<\alpha,\beta\>e^{\beta}_{m+n}$$ for
$h\in \frak{h},\; \beta\in L,\; m,n\in \Z$. As $W$ is irreducible,
$e^{\alpha}_{k}w$ generates $W$ by operators $h(m),\; e^{\beta}_{m}$
for $h\in \frak{h},\; \beta\in L,\; m\in \Z$. Then it follows that
$\alpha(0)$ acts semisimply on $W$ with only integer eigenvalues and
that $\alpha(n)$ for $n\ge 1$ act locally nilpotently.
\end{proof}

Let $\{\alpha_{1},\dots,\alpha_{r}\}\subset L$ be such that
$\{\alpha_{1},\dots,\alpha_{r}\}$ is a basis for the subspace $\C L$
of $\frak{h}$.  Then extend $\{\alpha_{1},\dots,\alpha_{r}\}$ to a
basis $\{ \alpha_{1},\dots,\alpha_{r},u_{1},\dots,u_{s}\}$ of
$\frak{h}$. Let $\{\beta_{1},\dots,\beta_{r},v_{1},\dots,v_{s}\}$ be
the dual basis. Set
\begin{eqnarray}
L^{o}=\{\lambda\in \frak{h}\;|\; \<\alpha,\lambda\>\in \Z\ \ \
\mbox{ for all }\alpha\in L \}.
\end{eqnarray}
Consider subalgebras of $\hat{\frak{h}}$:
\begin{eqnarray*}
{\mathcal{L}}_{0}&=&\sum_{i=1}^{r}\sum_{n\ge 1}(\C \alpha_{i}(n)+\C
\beta_{i}(-n))+\C {\bf k},\\
{\mathcal{L}}_{1}&=&\sum_{j=1}^{s}\sum_{n\ge 1}(\C u_{j}(n)+\C
v_{j}(-n))+\C {\bf k}.
\end{eqnarray*}
Then $\hat{\frak{h}}=({\mathcal{L}}_{0}+{\mathcal{L}}_{1})\oplus
\frak{h}(0)$ with ${\mathcal{L}}_{0}$, ${\mathcal{L}}_{1}$, and
$\frak{h}(0)$ mutually commuting. If $\C L=\frak{h}$, we have
${\mathcal{L}}_{1}=\C {\bf k}$ and if $L=0$, we have
${\mathcal{L}}_{0}=\C {\bf k}$. Otherwise, ${\mathcal{L}}_{0}$ and
${\mathcal{L}}_{1}$ are infinite-dimensional Heisenberg algebras.

Let $M_{0}(1)$ $(=\C[\beta_{i}(-n)\;| \; 1\le i\le d,\; n\ge 1]$ as
a vector space) denote the canonical irreducible
${\mathcal{L}}_{0}$-module of level $1$ on which $\alpha_{i}(n)$
acts locally nilpotent for $1\le i\le r,\; n\ge 1$. Let $U_{1}$ be
an irreducible  ${\mathcal{L}}_{1}$-module of level $1$ which is
restricted in the sense that for every $w\in U_{1}$ and $1\le j\le
s$, $u_{j}(n)w=0$ for $n$ sufficiently large.  Let $\lambda\in
L^{o}$ and denote by $\C e^{\lambda}$ (where $e^{\lambda}$ is just a
symbol) the $1$-dimensional $\frak{h}(0)$-module with $h(0)$ acting
as scalar $\<\lambda,h\>$ for $h\in \frak{h}$.
 Set
$$M(1,\lambda,
U_{1})=\C e^{\lambda}\otimes U_{1}\otimes M_{0}(1).$$ Then
$M(1,\lambda, U_{1})$ is an $\hat{\frak{h}}$-module of level $1$
with $\frak{h}(0)$ acting on $\C e^{\lambda}$, ${\mathcal{L}}_{1}$
acting on $U_{1}$, and ${\mathcal{L}}_{0}$ acting on $M_{0}(1)$.

\bp{plast} The defined $\hat{\frak{h}}$-module $M(1,\lambda,U_{1})$
satisfies Condition $C_{L}$ and is irreducible. On the other hand,
every irreducible $\hat{\frak{h}}$-module satisfying Condition
$C_{L}$ is isomorphic to $M(1,\lambda, U_{1})$ for some $\lambda\in
L^{o}$ and some irreducible restricted ${\mathcal{L}}_{1}$-module
$U_{1}$ of level $1$. \ep

\begin{proof} It is evident that $M(1,\lambda,U_{1})$
satisfies Condition $C_{L}$. Notice that the Schur lemma holds for
the ${\mathcal{L}}_{0}$-module $M_{0}(1)$ as $M_{0}(1)$ is of
countable dimension over $\C$. Then it follows that
$M(1,\lambda,U_{1})$ is an irreducible $\hat{\frak{h}}$-module. Let
$U$ be an irreducible $\hat{\frak{h}}$-module satisfying Condition
$C_{L}$. By Theorem \ref{theisenberg}, $U$ viewed as an
${\mathcal{L}}_{0}$-module is completely reducible with each
irreducible submodule isomorphic to $M_{0}(1)$. Consequently,
$U=M_{0}(1)\otimes U^{1},$ where
$U^{1}=\Hom_{\mathcal{L}_{0}}(M_{0}(1),U)$ is naturally an
$({\mathcal{L}}_{1}+\frak{h}(0))$-module of level $1$. With $U$ an
irreducible restricted $\hat{\frak{h}}$-module, $U^{1}$ irreducible
and restricted. As $\frak{h}(0)$ is central in $\hat{\frak{h}}$ and
commutes with ${\mathcal{L}}_{1}$, there exists $\lambda\in
\frak{h}$ such that $h(0)$ acts as scalar $\<\lambda,h\>$  on
$U^{1}$ for $h\in \frak{h}$, and $U^{1}$ is an irreducible
${\mathcal{L}}_{1}$-module of level $1$. From Condition $C_{L}$, we
have $\<\lambda,\alpha_{i}\>\in \Z$ for $1\le i\le r$. Thus
$\lambda\in L^{o}.$ Taking $U_{1}=U^{1}$ viewed as an
${\mathcal{L}}_{1}$-module, we have $U\simeq M(1,\lambda,U_{1})$.
\end{proof}

Let $\lambda\in L^{o}$ and $U_{1}$ an irreducible restricted
${\mathcal{L}}_{1}$-module of level $1$. In view of Proposition
\ref{plast} and Theorem \ref{tmain}, we have an irreducible
$V_{(\frak{h},L)}$-module $V_{(\frak{h},L)}(U)$ with
$U=M(1,\lambda,U_{1})$. We denote this module by
$V_{(\frak{h},L)}(\lambda,U_{1})$, where
$$V_{(\frak{h},L)}(\lambda,U_{1})=\C e^{\lambda}\otimes \C[L]\otimes
U_{1}\otimes M_{0}(1)$$ as a vector space.  Set
$${\frak{u}}(0)=\C u_{1}(0)\oplus \cdots \oplus \C u_{s}(0)\subset
\frak{h}(0)\subset \hat{\frak{h}}.$$

\bt{tconverse} Let $W$ be an irreducible $V_{(\frak{h},L)}$-module
on which $\frak{u}(0)$ acts semisimply. Then $W$ is isomorphic to
$V_{(\frak{h},L)}(\lambda,U_{1})$ for some $\lambda\in L^{o}$ and
for some irreducible restricted ${\mathcal{L}}_{1}$-module $U_{1}$
of level $1$. \et

\begin{proof} For $\lambda\in \frak{h}$, set
$$W_{\lambda}=\{ w\in W\;|\; h(0)w=\<\lambda,h\>w\ \ \mbox{ for
}h\in \frak{h}\}.$$ Because $\alpha_{1}(0),\dots,\alpha_{r}(0)$ are
semisimple by Lemma \ref{lconverse} and $\frak{u}(0)$ is assumed to
be semisimple, we have $W=\coprod_{\lambda\in \frak{h}}W_{\lambda}$.
As $[\frak{h}(0),\hat{\frak{h}}]=0$, $W_{\lambda}$ are
$\hat{\frak{h}}$-submodules of $W$. For $\alpha\in L,\; h\in
\frak{h},\; v\in M(1)$, we have
$$[h(0),Y_{W}(e^{\alpha}\otimes v,x)]=Y_{W}(h(0)(e^{\alpha}\otimes
v),x) =\<h,\alpha\>Y_{W}(e^{\alpha}\otimes v,x).$$ Consequently,
\begin{eqnarray}\label{especial-bracket}
Y_{W}(e^{\alpha}\otimes v,x)W_{\lambda}\subset
W_{\lambda+\alpha}[[x,x^{-1}]]
\end{eqnarray}
for $\alpha\in L,\; v\in M(1),\; \lambda\in \frak{h}$. Suppose that
$W_{\lambda_{0}}\ne 0$ and let $0\ne w\in W_{\lambda_{0}}$. By a
result of \cite{dm} and \cite{li-thesis}, the linear span of $\{
a_{m}w\;|\; a\in V_{(\frak{h},L)},\; m\in \Z\}$ is a
$V_{(\frak{h},L)}$-submodule of $W$. Consequently,
$$W={\rm span}\{ a_{m}w\;|\; a\in V_{(\frak{h},L)},\; m\in \Z\}.$$
Combining this with (\ref{especial-bracket}) and the decomposition
$W=\coprod_{\lambda\in \frak{h}}W_{\lambda}$, we get
$$W_{\lambda_{0}}={\rm span}\{ v_{m}w\;|\; v\in M(1),\; m\in \Z\}
\ \ \mbox{ and }\ \ W=\coprod_{\alpha\in L}W_{\lambda_{0}+\alpha}.$$
It follows that $W_{\lambda_{0}}$ is an irreducible module for
$M(1)$ viewed as a vertex algebra. As $\frak{h}$ generates $M(1)$ as
a vertex algebra, $W_{\lambda_{0}}$ is an irreducible
$\hat{\frak{h}}$-module. In view of Lemma \ref{lconverse},
$W_{\lambda_{0}}$ satisfies Condition $C_{L}$. We are going to prove
that $W\simeq V_{(\frak{h},L)}(W_{\lambda_{0}})$.

Let $\alpha\in L$. We have
\begin{eqnarray*}
& &\frac{d}{dx}Y_{W}(e^{\alpha},x)=Y_{W}(L(-1)e^{\alpha},x)=
Y_{W}(\alpha(-1)e^{\alpha},x)\\
& &\ \ \ \ \
=\alpha(x)^{+}Y_{W}(e^{\alpha},x)+Y_{W}(e^{\alpha},x)\alpha(x)^{-},
\end{eqnarray*}
where $\alpha(x)^{+}=\sum_{n<0}\alpha(n)x^{-n-1}$ and
$\alpha(x)^{-}=\sum_{n\ge 0}\alpha(n)x^{-n-1}$. We also have
$\frac{d}{dx}E^{-}(\alpha,x)=-E^{-}(\alpha,x)\alpha(x)^{+}$, and
$\frac{d}{dx}E^{+}(\alpha,x)x^{-\alpha(0)}=-\alpha(x)^{-}E^{+}(\alpha,x)x^{-\alpha(0)}$.
Using these relations we obtain
$$\frac{d}{dx}E^{-}(\alpha,x)Y_{W}(e^{\alpha},x)E^{+}(\alpha,x)x^{-\alpha(0)}=0.$$
Then we set
$$E_{\alpha}=E^{-}(\alpha,x)Y_{W}(e^{\alpha},x)E^{+}(\alpha,x)x^{-\alpha(0)}\in \End _{\C}W.$$
Define a linear map $\psi:
V_{(\frak{h},L)}(W_{\lambda_{0}})\rightarrow W$ by
\begin{eqnarray}
\psi( e^{\alpha}\otimes u)=E_{\alpha}u\ \ \ \mbox{ for }\alpha\in
L,\; u\in W_{\lambda_{0}}.
\end{eqnarray}
We are going to prove that $\psi$ is a $V_{(\frak{h},L)}$-module
isomorphism.

First, we establish some properties for $E_{\alpha}$ with $\alpha\in
L$. For $h\in \frak{h}$, we have
\begin{eqnarray}
[h(0),E_{\alpha}]=E^{-}(\alpha,x)[h(0),Y_{W}(e^{\alpha},x)]E^{+}(\alpha,x)x^{-\alpha(0)}
=\<h,\alpha\>E_{\alpha}.
\end{eqnarray}
As $h(i)e^{\alpha}=\delta_{i,0}\<h,\alpha\>e^{\alpha}$ for $i\ge 0$,
we have
$$[h(n),Y_{W}(e^{\alpha},x)]=\<h,\alpha\>x^{n}Y_{W}(e^{\alpha},x)
\ \ \ \ \mbox{ for }n\in \Z.$$ If $n>0$, from \cite{flm}
(Proposition 4.1.1) we also have
\begin{eqnarray*}
& &[h(n),E^{-}(\alpha,x)]=-\<h,\alpha\>x^{n}E^{-}(\alpha,x),\\
& &[h(n),E^{+}(\alpha,x)]=0=[h(n),\alpha(0)].
\end{eqnarray*}
Then we get $[h(n),E_{\alpha}]=0$. Similarly, we have $[E_{\alpha},
h(n)]=0$ for $h\in \frak{h},\; n< 0$.

{}From the definition of $E_{\alpha}$ we have
\begin{eqnarray}\label{eYW}
Y_{W}(e^{\alpha},x)=E^{-}(-\alpha,x)E^{+}(-\alpha,x)E_{\alpha}x^{\alpha(0)}.
\end{eqnarray}
For $\alpha,\beta\in L$, from (\ref{eY-formula-special}) we have
$$e^{\alpha}_{-\<\alpha,\beta\>-1}e^{\beta}=\varepsilon(\alpha,\beta)e^{\alpha+\beta}
\ \ \mbox{ and }\ \ e^{\alpha}_{m}e^{\beta}=0\ \ \mbox{ for }m\ge
-\<\alpha,\beta\>.$$ Let $0\ne w\in W_{\lambda_{0}}$ be such that
$\gamma(n)w=0$ for $\gamma\in L,\; n\ge 1$. Using (\ref{eYW}) we
have
$$e^{\gamma}_{-\<\gamma,\lambda_{0}\>-1}w=E_{\gamma}w \ \ \mbox{ and }\ \ e^{\gamma}_{m}w=0
\ \ \mbox{ for }\gamma\in L, \; m\ge -\<\gamma,\lambda_{0}\>.$$ In
particular, this is true for $\gamma=\alpha$, or $\beta$. Combining
this with Jacobi identity we get
$$(x_{0}+x_{2})^{-\<\alpha,\lambda_{0}\>}Y_{W}(e^{\alpha},x_{0}+x_{2})Y_{W}(e^{\beta},x_{2})w
=(x_{2}+x_{0})^{-\<\alpha,\lambda_{0}\>}Y_{W}(Y(e^{\alpha},x_{0})e^{\beta},x_{2})w.$$
Then by applying
$\Res_{x_{0}}\Res_{x_{2}}x_{0}^{-\<\alpha,\beta\>-1}x_{2}^{-\<\beta,\lambda_{0}\>-1}$
we obtain
\begin{eqnarray}
E_{\alpha}E_{\beta}=\varepsilon(\alpha,\beta)E_{\alpha+\beta}.
\end{eqnarray}

Finally, we are ready to finish the proof.  We have
$$\psi (h(0)(e^{\alpha}\otimes u))
=\<\alpha,h\>E_{\alpha}u+E_{\alpha}h(0)u=h(0)E_{\alpha}u
=h(0)\psi(e^{\alpha}\otimes u).$$
Then $\psi$ is an
$\hat{\frak{h}}$-module homomorphism. Furthermore, we have
\begin{eqnarray*}
\psi\left( Y_{W}(e^{\alpha},x)(e^{\beta}\otimes u)\right)
&=&\psi\left(\varepsilon(\alpha,\beta)e^{\alpha+\beta}\otimes
x^{\<\alpha,\beta\>}E^{-}(-\alpha,x)E^{+}(-\alpha,x)x^{\alpha(0)}u\right)\\
&=&\varepsilon(\alpha,\beta)x^{\<\alpha,\beta\>}
E_{\alpha+\beta}E^{-}(-\alpha,x)E^{+}(-\alpha,x)x^{\alpha(0)}u\\
&=&x^{\<\alpha,\beta\>}E_{\alpha}E_{\beta}E^{-}(-\alpha,x)E^{+}(-\alpha,x)x^{\alpha(0)}u\\
&=&E^{-}(-\alpha,x)E^{+}(-\alpha,x)E_{\alpha}
x^{\alpha(0)}E_{\beta}u\\
&=&Y_{W}(e^{\alpha},x)\psi(e^{\beta}\otimes u).
\end{eqnarray*}
As $\frak{h}$ and $e^{\alpha}\ (\alpha\in L)$ generate
$V_{(\frak{h},L)}$ as a vertex algebra, $\psi$ is a
$V_{(\frak{h},L)}$-module isomorphism.
\end{proof}

\bp{pequivalence} Let $\lambda_{1},\lambda_{2}\in L^{o}$ and let
$U_{1}$ and $U_{2}$ be irreducible restricted
${\mathcal{L}}_{1}$-modules of level $1$. Then
$V_{(\frak{h},L)}(\lambda_{1},U_{1})\simeq
V_{(\frak{h},L)}(\lambda_{2},U_{2})$ if and only if
$\lambda_{1}+L=\lambda_{2}+L$ and $U_{1}\simeq U_{2}$. \ep

\begin{proof} Note that for each $\alpha\in L$,
$\C e^{\alpha} \otimes \C e^{\lambda_{i}}\otimes U_{i}\otimes
M_{0}(1)$ is an $\hat{\frak{h}}$-submodule of
$V_{(\frak{h},L)}(\lambda_{i},U_{i})$ and
$$\C e^{\alpha} \otimes \C e^{\lambda_{i}}\otimes U_{i}\otimes
M_{0}(1)\simeq \C e^{\lambda_{i}+\alpha}\otimes U_{i}\otimes
M_{0}(1)=M(1,\lambda_{i}+\alpha,U_{i}).$$ We see that the set of
$\frak{h}(0)$-weights of $V_{(\frak{h},L)}(\lambda_{i},U_{i})$ is
$\lambda_{i}+ L$. If $V_{(\frak{h},L)}(\lambda_{1},U_{1})\simeq
V_{(\frak{h},L)}(\lambda_{2},U_{2})$, we must have
$\lambda_{1}+L=\lambda_{2}+L$ and $U_{1}\otimes M_{0}(1)\simeq
U_{2}\otimes M_{0}(1)$ as
$({\mathcal{L}}_{0}+{\mathcal{L}}_{1})$-modules, which implies that
$U_{1}\simeq U_{2}$.

On the other hand, assume $\lambda_{1}+L=\lambda_{2}+L$ and
$U_{1}=U_{2}$. Then $\lambda_{1}=\lambda_{2}+\gamma$ for some
$\gamma\in L$. Define a linear isomorphism $\theta:
V_{(\frak{h},L)}(\lambda_{1},U_{1})\rightarrow
V_{(\frak{h},L)}(\lambda_{2},U_{2})$ by
$$\theta (e^{\beta}\otimes e^{\lambda_{1}}\otimes w)
=\varepsilon(\beta,\gamma) (e^{\beta+\gamma}\otimes
e^{\lambda_{2}}\otimes w)$$ for $\beta\in L,\; w\in U_{1}\otimes
M_{0}(1)$. It is clear that $\theta$ is an $\hat{\frak{h}}$-module
isomorphism. Furthermore, for $\alpha\in L$, we have
\begin{eqnarray*}
& &\theta( Y_{W}(e^{\alpha},x)(e^{\beta}\otimes
e^{\lambda_{1}}\otimes w))\\
&=&\theta\left( \varepsilon(\alpha,\beta)e^{\alpha+\beta} \otimes
x^{\<\alpha,\beta+\lambda_{1}\>}e^{\lambda_{1}}
\otimes E^{-}(-\alpha,x)E^{+}(-\alpha,x)w \right)\\
&=&\varepsilon(\alpha,\beta)\varepsilon(\alpha+\beta,\gamma)e^{\alpha+\beta+\gamma}
\otimes x^{\<\alpha,\beta+\lambda_{1}\>}e^{\lambda_{2}}
\otimes E^{-}(-\alpha,x)E^{+}(-\alpha,x)w\\
&=&\varepsilon(\alpha,\beta+\gamma)\varepsilon(\beta,\gamma)e^{\alpha+\beta+\gamma}
\otimes x^{\<\alpha,\beta+\gamma+\lambda_{2}\>}e^{\lambda_{2}}
\otimes E^{-}(-\alpha,x)E^{+}(-\alpha,x)w\\
&=&Y_{W}(e^{\alpha},x)\theta(e^{\beta}\otimes e^{\lambda_{1}}\otimes
w).
\end{eqnarray*}
Since $\frak{h}$ and
$e^{\alpha}$ for $\alpha\in L$ generate $V_{(\frak{h},L)}$ as a
vertex algebra, $\theta$ is a $V_{(\frak{h},L)}$-module
homomorphism.
\end{proof}

\br{rnondegenerate-lattice} {\em Assume that $L$ is a nondegenerate
even lattice with $\frak{h}=\C\otimes_{\Z}L$. Let $U$ be an
$\hat{\frak{h}}$-module satisfying Condition $C_{L}$. Then for any
$h\in \frak{h},\; n\ge 1$, $h(n)$ is locally nilpotent and $h(0)$ is
semisimple on $U$. It follows that $U$ is completely reducible.
Since $\alpha(0)$ for $\alpha\in L$ are assumed to have only integer
eigenvalues, $U$ is a direct sum of irreducible modules isomorphic
to $M(1)\otimes \C e^{\beta}$ for $\beta\in L^{o}$.
 If
$U$ is irreducible, we have $U=M(1)\otimes \C e^{\beta}$ for some
$\beta\in L^{o}$. Thus Theorem \ref{tmain} recovers the
corresponding results of \cite{flm} with a different proof. In this
case, the assumption in Theorem \ref{tconverse} that $\frak{u}(0)$
is semisimple is vacuous. Then we recover Dong's corresponding
result in \cite{dong1} (cf. \cite{dlm-reg}, \cite{ll}). }\er

\br{rdegenerate-lattice} {\em  Let $(\frak{h},L)$ be a pair as in
Theorems \ref{tmain} and \ref{tconverse} such that $\C L\ne
\frak{h}$. In this case, ${\mathcal{L}}_{1}$ is infinite-dimensional
and from Section 2 there are many irreducible restricted
${\mathcal{L}}_{1}$-modules of level $1$ besides the one on the
polynomial algebra, which give rise to families of irreducible
$V_{(\frak{h},L)}$-modules.} \er

\section{A characterization of vertex algebras $V_{(\frak{h},L)}$
with $\<\cdot,\cdot\>|_{L}=0$}

In this section we study vertex algebras $V_{(\frak{h},L)}$ with
$\<\cdot,\cdot\>|_{L}=0$. In this case, we give a characterization
of the vertex algebras  in terms of a certain affine Lie algebra.

Let $\frak{h}$ be a finite-dimensional vector space equipped with a
nondegenerate symmetric bilinear form $\<\cdot,\cdot\>$ and let
$L\subset \frak{h}$ a free abelian group with $\<\alpha,\beta\>=0$
for $\alpha,\beta\in L$. View $\frak{h}$ and the group algebra
$\C[L]$ as abelian Lie algebras. Let $\frak{h}$ act on $\C[L]$ by
$$h\cdot e^{\alpha}=\<h,\alpha\>e^{\alpha}\ \ \ \mbox{ for }h\in
\frak{h},\; \alpha\in L.$$ Then $\frak{h}$ acts on $\C[L]$ (viewed
as a Lie algebra) by derivations. Form the cross product Lie algebra
$\frak{p}=\frak{h}\ltimes \C[L]$ and extend the bilinear form on
$\frak{h}$ to $\frak{p}$ by
\begin{eqnarray}
\<h+u,h'+v\>=\<h,h'\> \ \ \ \ \mbox{ for }h,h'\in \frak{h},\; u,v\in
\C[L].
\end{eqnarray}
This form is symmetric and invariant. Then we have an affine Lie
algebra
$$\hat{\frak{p}}=\frak{p}\otimes \C[t,t^{-1}]\oplus \C {\bf k},$$
where ${\bf k}$ is central and
\begin{eqnarray*}
& &[h(m),h'(n)]=m\<h,h'\>\delta_{m+n,0}{\bf k},\\
& &[h(m),e^{\alpha}(n)]=\<h,\alpha\>e^{\alpha}(m+n),\\
& &[e^{\alpha}(m),e^{\beta}(n)]=0
\end{eqnarray*}
for $h,h'\in \frak{h},\; \alpha,\beta\in L,\; m,n\in \Z$. Let $\ell$
be a complex number. Denote by $\C_{\ell}$ the $1$-dimensional
$(\frak{p}\otimes \C[t]+\C {\bf k})$-module $\C$ with
$\frak{p}\otimes \C[t]$ acting trivially and with ${\bf k}$ acting
as scalar $\ell$. Form the induced $\hat{\frak{p}}$-module
$$V_{\hat{\frak{p}}}(\ell,0)
=U(\hat{\frak{p}})\otimes_{U(\frak{p}\otimes \C[t]+\C {\bf k})}
\C_{\ell}.$$ We identify $\frak{p}$ as a subspace of
$V_{\hat{\frak{p}}}(\ell,0)$ through the linear map $a\mapsto
a(-1){\bf 1}$, where ${\bf 1}=1\otimes 1$. In particular, we
identify $\alpha$ with $\alpha(-1){\bf 1}$ and $e^{\alpha}$ with
$e^{\alpha}(-1){\bf 1}$ for $\alpha\in L$. Then there exists a
vertex algebra structure on $V_{\hat{\frak{p}}}(\ell,0)$ with ${\bf
1}$ as the vacuum vector and with $Y(a,x)=a(x)$ for $a\in \frak{p}$
(cf. \cite{ll}). Recall that $\D$ is the linear operator on
$V_{\hat{\frak{p}}}(\ell,0)$ defined by $\D v=v_{-2}{\bf 1}$. Note
that the normalized $2$-cocycle $\varepsilon: L\times L\rightarrow
\C^{\times}$ defined in Remark \ref{rcocycle-exist} is trivial in
the sense that $\varepsilon(\alpha,\beta)=1$ for all
$\alpha,\beta\in L$.

 \bt{tchara} The vertex algebra $V_{(\frak{h},L)}$
is isomorphic to the quotient vertex algebra of
$V_{\hat{\frak{p}}}(1,0)$ modulo the ideal $J$ generated by the
elements
$$ e^{0}-{\bf 1},\ \ \ e^{\alpha}_{-1}e^{\beta}-e^{\alpha+\beta},
\ \ \ \D e^{\alpha}-\alpha(-1)e^{\alpha}\ \ \ \mbox{ for
}\alpha,\beta\in L.$$ \et

\begin{proof} First, we show that $V_{(\frak{h},L)}$ is a
$\hat{\frak{p}}$-module of level $1$ with $h(x)=Y(h,x)$ and
$e^{\alpha}(x)=Y(e^{\alpha},x)$ for $h\in \frak{h},\; \alpha \in L$.
We know that $V_{(\frak{h},L)}$ is an $\hat{\frak{h}}$-module of
level $1$ with $h(x)=Y(h,x)$. For $\alpha,\beta\in L$, since
$\<\alpha,\beta\>=0$, we have
$$Y(e^{\alpha},x)e^{\beta}
=E^{-}(-\alpha,x)E^{+}(-\alpha,x)e^{\alpha}x^{\alpha(0)}\cdot
e^{\beta} =E^{-}(-\alpha,x)e^{\alpha+\beta},$$ which contains only
nonnegative powers of $x$, so that
$$[Y(e^{\alpha},x_{1}),Y(e^{\beta},x_{2})]=0.$$
For $h\in \frak{h}$, since
$h(n)e^{\alpha}=\delta_{n,0}\<h,\alpha\>e^{\alpha}$ for $n\ge 0$, we
have
$$[Y(h,x_{1}),Y(e^{\alpha},x_{2})]
=\<h,\alpha\>x_{2}^{-1}\delta\left(\frac{x_{1}}{x_{2}}\right)Y(e^{\alpha},x_{2}).$$
Then $V_{(\frak{h},L)}$ is a $\hat{\frak{p}}$-module of level $1$.
Clearly, ${\bf 1}$ generates $V_{(\frak{h},L)}$ as a
$\hat{\frak{p}}$-module with $a(n){\bf 1}=0$ for $a\in \frak{p},\;
n\ge 0$. It follows that there exists a unique
$\hat{\frak{p}}$-module homomorphism $\psi$ from
$V_{\hat{\frak{p}}}(1,0)$ onto $V_{(\frak{h},L)}$, sending the
vacuum vector to the vacuum vector. That is, $\psi({\bf 1})={\bf 1}$
and
$$\psi(a_{m}v)=a_{m}\psi(v)\ \ \ \mbox{ for }a\in \frak{p},
\ v\in V_{\hat{\frak{p}}}(1,0),\ m\in \Z.$$ As $\frak{p}$ generates
$V_{\hat{\frak{p}}}(1,0)$ as a vertex algebra, it follows that
$\psi$ is a vertex algebra homomorphism. It is easy to see that the
following relations hold in $V_{(\frak{h},L)}$:
$$e^{0}={\bf 1},\ \ e^{\alpha}_{-1}e^{\beta}=e^{\alpha+\beta},\ \
\D e^{\alpha}=e^{\alpha}_{-2}{\bf 1}=\alpha(-1)e^{\alpha}\ \ \
\mbox{ for }\alpha\in L,$$ so that the kernel of $\psi$ contains the
ideal $J$. Set $V=V_{\hat{\frak{p}}}(1,0)/J$ (the quotient vertex
algebra). Due to the linear map $\psi$, the linear map $a\in
\frak{p}\mapsto a(-1){\bf 1}\in V$ is also injective, so that
$\frak{p}$ can be identified as a subspace of $V$. Set
$$K=\sum_{\alpha\in L}U(\hat{\frak{h}})e^{\alpha}\subset V.$$
We shall show that $K=V$ by proving that $K$ is a vertex subalgebra
of $V$, containing all the generators. As $[\D,h(m)]=-mh(m-1)$ for
$h\in \frak{h},\; m\in \Z$ and $\D e^{\alpha}=\alpha(-1)e^{\alpha}$
for $\alpha\in L$, we see that $\D K\subset K$. For $\alpha,\beta\in
L$, we have $e^{\alpha}_{n}e^{\beta}=0$ for $n\ge 0$ and
$e^{\alpha}_{-1}e^{\beta}=e^{\alpha+\beta}$. For $n\ge 1$, we have
\begin{eqnarray*}
& &ne^{\alpha}_{-n-1}e^{\beta}=[\D, e^{\alpha}_{-n}]e^{\beta} =\D
e^{\alpha}_{-n}e^{\beta}-e^{\alpha}_{-n}\D e^{\beta} =\D
e^{\alpha}_{-n}e^{\beta}-e^{\alpha}_{-n}\beta(-1)e^{\beta}\\
& &\ \ \ \ =(\D-\beta(-1))e^{\alpha}_{-n}e^{\beta},
\end{eqnarray*}
noting that $[\beta(-1),e^{\alpha}_{m}]=0$ for $m\in \Z$. It follows
from induction that $e^{\alpha}_{-n-1}e^{\beta}\in K$ for $n\ge 0$.
Thus
$$e^{\alpha}_{m}e^{\beta}\in K\ \ \ \mbox{ for }\alpha,\beta,\; m\in
\Z.$$ For $h\in \frak{h},\;\alpha\in L,\; \;m, n\in \Z$, we have
$$e^{\alpha}_{m}h(n)=h(n)e^{\alpha}_{m}+\<\alpha,h\>e^{\alpha}_{m+n}.$$
Again, it follows from induction that $e^{\alpha}_{m}K\subset K$ for
$m\in \Z$. Thus $K$ is a $\hat{\frak{p}}$-submodule of $V$,
containing the vacuum vector ${\bf 1}$. Consequently, $K=V$. For
$\alpha\in L$, $U(\hat{\frak{h}})e^{\alpha}$ is an irreducible
$\hat{\frak{h}}$-module of level $1$, which is isomorphic to
$M(1)\otimes \C e^{\alpha}$. Then it follows that $\psi$ is an
isomorphism.
\end{proof}

\bp{pmodue} For any $V_{(\frak{h},L)}$-module $(W,Y_{W})$, $W$ is a
restricted $\hat{\frak{p}}$-module of level $1$ with
$a(x)=Y_{W}(a,x)$ for $a\in \frak{p}$, satisfying the condition that
\begin{eqnarray}
e^{0}(x)=1,\ \ \ e^{\alpha+\beta}(x)=e^{\alpha}(x)e^{\beta}(x),\ \ \
\frac{d}{dx}e^{\alpha}(x)=\alpha(x)e^{\alpha}(x)
\end{eqnarray}
for $\alpha\in L$. On the other hand, if $W$ is a restricted
$\hat{\frak{p}}$-module of level $1$, satisfying the above
condition, then there exists a $V_{(\frak{h},L)}$-module structure
$Y_{W}$ such that $Y_{W}(a,x)=a(x)$ for $a\in \frak{p}$.
 \ep

\begin{proof} Let $(W,Y_{W})$ be a $V_{(\frak{h},L)}$-module.
Then $W$ is naturally a $V_{\hat{\frak{p}}}(1,0)$-module. We have
$$e^{0}(x)=Y_{W}(e^{0},x)=Y_{W}({\bf 1},x)=1,$$
\begin{eqnarray*}
& &\frac{d}{dx}e^{\alpha}(x)=\frac{d}{dx}Y_{W}(e^{\alpha},x)
=Y_{W}(\D e^{\alpha},x)=Y_{W}(\alpha(-1)e^{\alpha},x)
=Y_{W}(\alpha,x)Y_{W}(e^{\alpha},x)
\ \ \ \ \\
& &\ \ \ \ \ \ \ \ =\alpha(x)e^{\alpha}(x),\\
& &e^{\alpha}(x)e^{\beta}(x)=Y_{W}(e^{\alpha},x)Y_{W}(e^{\beta},x)
=Y_{W}(e^{\alpha}_{-1}e^{\beta},x)=Y_{W}(e^{\alpha+\beta},x)=e^{\alpha+\beta}(x).
\end{eqnarray*}

Conversely, let $W$ be a restricted $\hat{\frak{p}}$-module of level
$1$, satisfying the conditions. Then $W$ is a
$V_{\hat{\frak{p}}}(1,0)$-module with $Y_{W}(a,x)=a(x)$ for $a\in
\frak{p}$. We have
\begin{eqnarray*}
Y_{W}(e^{0}-{\bf 1},x)&=&e^{0}(x)-1=0,\\
Y_{W}(e^{\alpha}_{-1}e^{\beta}-e^{\alpha+\beta},x)
&=&Y_{W}(e^{\alpha},x)Y_{W}(e^{\beta},x)-Y_{W}(e^{\alpha+\beta},x)\\
&=&e^{\alpha}(x)e^{\beta}(x)-e^{\alpha+\beta}(x)=0,\\
Y_{W}(\D e^{\alpha}-\alpha(-1)e^{\alpha},x)
&=&\frac{d}{dx}Y_{W}(e^{\alpha},x)-Y_{W}(\alpha,x)Y_{W}(e^{\alpha},x)\\
&=&\frac{d}{dx}e^{\alpha}(x)-\alpha(x)e^{\alpha}(x)=0.
\end{eqnarray*}
Then $W$ is naturally a module for the quotient vertex algebra
$V_{\hat{\frak{p}}}(1,0)/J$, where $J$ is the ideal generated by the
vectors
$$e^{0}-{\bf 1},\ \ e^{\alpha}_{-1}e^{\beta}-e^{\alpha+\beta},\;\;
\D e^{\alpha}-\alpha(-1)e^{\alpha}$$ for $\alpha,\beta\in L$. In
view of Theorem \ref{tchara}, $W$ is a $V_{(\frak{h},L)}$-module.
\end{proof}

\end{document}